\documentclass[10pt,letter]{amsart}
\usepackage{latexsym}
\usepackage{amscd, amsfonts, eucal, mathrsfs, amsmath, amssymb, amsthm, stmaryrd}
\usepackage{fontenc}%[OT2,OT1]{fontenc}
\usepackage{newcent}

\input xy
\xyoption{all}

\usepackage{fullpage}

\newcommand{\field}[1]{\mathbf #1}
\newcommand{\mf}[1]{\mathfrak #1}

\newcommand{\ms}[1]{\mathscr #1}
\newcommand{\widebar}[1]{\overline{#1}}

\newcommand{\R}{\field R}
\newcommand{\LLL}{\field L}

\newcommand{\F}{\field F}
\newcommand{\Z}{\field Z}
\newcommand{\Q}{\field Q}
\newcommand{\N}{\field N}
\newcommand{\V}{\field V}

\newcommand{\simto}{\stackrel{\sim}{\to}}
\newcommand{\eps}{\varepsilon}
\renewcommand{\phi}{\varphi}

\renewcommand{\hom}{\operatorname{Hom}}
\newcommand{\Hom}{\operatorname{Hom}}
\newcommand{\shom}{\ms H\!om}

\DeclareMathOperator{\uhom}{\underline{Hom}}
\newcommand{\rshom}{\mathbf{R}\shom}

\DeclareMathOperator{\rhom}{\operatorname{{\bf R}Hom}}
\newcommand{\send}{\ms E\!nd}

\newcommand{\spec}{\operatorname{Spec}}

\newcommand{\A}{\field A}

\DeclareMathOperator{\Pic}{Pic}

\newcommand{\td}{\operatorname{Td}}
\DeclareMathOperator{\chern}{ch}

\DeclareMathOperator{\supp}{Supp}
\DeclareMathOperator{\pr}{pr}

\DeclareMathOperator{\Quot}{Quot}

\newcommand{\m}{\boldsymbol{\mu}}

\newcommand{\G}{\field G} %for the multiplicative and additive groups

\renewcommand{\H}{\operatorname{H}}

\DeclareMathOperator{\ext}{\operatorname{Ext}}

\DeclareMathOperator{\Coh}{\operatorname{Coh}}

\DeclareMathOperator{\D}{\operatorname{\bf D}}

\renewcommand{\]}{]\!\hspace{0.03em}]}

\DeclareMathOperator*{\tensor}{\otimes}
\DeclareMathOperator*{\ltensor}{\stackrel{\field L}{\otimes}}

\renewcommand{\N}{\operatorname{N}}

\DeclareMathOperator{\rk}{\operatorname{rk}}

\newcommand{\surj}{\twoheadrightarrow}
\newcommand{\inj}{\hookrightarrow}

\DeclareMathOperator{\End}{\operatorname{End}}
\DeclareMathOperator{\aut}{\operatorname{Aut}}
\DeclareMathOperator{\Aut}{\operatorname{Aut}}

\DeclareMathOperator{\M}{\operatorname{M}}

\renewcommand{\sp}{\text{\rm sp}} %space locus: where a stack is a space
\DeclareMathOperator{\sing}{Sing}
\newcommand{\dom}{\text{\rm dom}}
\newcommand{\softboundary}{\partial^{\text{\rm soft}}}
\newcommand{\soft}{\text{\rm soft}}

\newcommand{\Tw}{\mathbf{Tw}}
\DeclareMathOperator{\mTw}{Tw}

\DeclareMathOperator{\B}{\operatorname{\mathsf B\!}}

\DeclareMathOperator{\codim}{codim}
\DeclareMathOperator{\Filt}{Filt}

\newtheorem{lem}{Lemma}[subsubsection]

\makeatletter
\@addtoreset{lem}{section}
\@addtoreset{lem}{subsection}
\@addtoreset{lem}{subsubsection}
\@addtoreset{lem}{paragraph}
\makeatother

\renewcommand{\thelem}{\ifnum\value{subsubsection}>0{\thesubsubsection.\arabic{lem}}\else{\ifnum\value{subsection}>0{\thesubsection.\arabic{lem}}\else{\thesection.\arabic{lem}}\fi}\fi}

\newtheorem{thm}[lem]{Theorem}
\newtheorem*{theorem}{Theorem}

\newtheorem{prop}[lem]{Proposition}

\newtheorem{cor}[lem]{Corollary}

\newtheorem*{claim}{Claim}

\theoremstyle{definition}
\newtheorem{defn}[lem]{Definition}

\newtheorem{example}[lem]{Example}
\newtheorem{hyp}[lem]{Hypothesis}

\theoremstyle{remark}
\newtheorem{remark}[lem]{Remark}

\newtheorem{notation}[lem]{Notation}

\numberwithin{equation}{lem}

\author{Max Lieblich}
\address{Fine Hall, Washington Road, Princeton NJ 08544-1000}
\email{lieblich@math.princeton.edu}
\title{Moduli of twisted orbifold sheaves}
\begin{document}
%\bibliographystyle{plain}

%\centerline{\it Dedicated to the memory of Joseph J.\ Katz}

\begin{abstract}
  We study stacks of slope-semistable twisted sheaves on orbisurfaces
  with projective coarse spaces and prove that in certain cases they
  have many of the asymptotic properties enjoyed by the moduli of
  slope-semistable sheaves on smooth projective surfaces.
\end{abstract}

\maketitle

\tableofcontents
\section{Introduction}
\label{sec:introduction}

The purpose of this paper is to streamline and broaden the moduli
theory of twisted sheaves on surfaces.  We also lay some foundations
for the study of moduli of (twisted) sheaves on stacky surfaces. In
particular, we will prove the following theorem.  Let $k$ be a field
and $r$ an integer which is invertible in $k$.

\begin{theorem}
  Let $\ms X\to X$ be an $\m_N$-gerbe on an $A_{N-1}$-orbisurface with
  projective coarse moduli space.  Given an invertible sheaf
  $L\in\Pic(X)$ and a positive integer $\Delta$, the stack
  $\Tw^{ss}_{\ms X/k}(r,L,\Delta)$ of totally regular slope-semistable
  $\ms X$-twisted sheaves of rank $r$, determinant $L$, and
  discriminant $\Delta$ is of finite type over $k$ and geometrically
  normal.  In addition, 
  \begin{enumerate}
  \item if there exists a locally free $\ms X$-twisted sheaf of rank
    $r$, then for any positive integer $B$, there is $\Delta>B$ such
    that $\Tw^{ss}_{\ms X/k}(r,L,\Delta)$ contains at least one
    geometrically integral connected component;
  \item if $X$ is a smooth projective surface, then for all
    sufficiently large $\Delta$, the stack $\Tw^{ss}_{\ms
      X/k}(r,L,\Delta)$ is geometrically integral whenever it is
    non-empty, and it is non-empty for infinitely many values of
    $\Delta$.  Moreover, if it is non-empty for $\Delta$, then it is
    non-empty for $\Delta+r\ell$ for any $\ell>0$.
  \end{enumerate}
\end{theorem}
The condition that $X$ is an $A_{N-1}$-orbisurface means, roughly,
that $X$ has isolated stacky points, each with stabilizer $\m_N$.  (We
also permit the locus of stacky points to be empty, in which case $X$
is just a smooth projective surface.)  The
condition that a sheaf is totally regular means that (the classes in
$K$-theory of) its (derived) representations at the stacky points of
$X$ are multiples of the standard representation of $\m_N$.

In previous work \cite{twisted-moduli}, we proved the second part of
the above theorem under the additional assumption that the class of
$\ms X$ in the Brauer group of $X\tensor_k\widebar k$ had order $N$
and that $r=N$.  The more general results proven here bring the
structure theorems for moduli of twisted sheaves on smooth projective
surfaces into alignment with those known for classical sheaves.  The first
part of the theorem uses a method shown to the author by Johan de
Jong, which applies to a large range of moduli problems, notably in
recent work of de Jong and Starr on rationally connected varieties.
Here the asympotic regularity properties of the moduli stacks lead to
the existence of geometrically irreducible components, but we cannot
yet conclude that the entire moduli problem becomes asymptotically irreducible.

The approach we take owes a debt to Adrian Langer, who, in a series of
papers (\cite{langer}, \cite{MR2085175}, and
\cite{langer-castelnuovo}) extended O'Grady's beautiful methods \cite{o'grady}
to all characteristics, clarifying the role of the Segre loci.  It is
O'Grady's now standard outline, reproduced in \cite{h-l} and
\cite{langer-castelnuovo}, that we have followed.  In several cases,
Langer's proofs carry over almost verbatim.  When this is the case, we
refer to his proofs and only indicate the changes that must be made in
the present context.

There are three main ways in which our proof differs from the standard
outline.  The first is that
the delicate numerical estimates needed for bounding the
dimensions of the Segre loci are altered by the replacement of the
classical Riemann-Roch formula by To\"en's orbifold Riemann-Roch
formula.  We include computations of the correction terms, which, in
the case of totally regular sheaves, do not disturb the main
properties of the relevant estimates.

The second difference between our proof and the standard outline is in
the terminal steps of the proof of irreducibility, where Langer's
effective computations \cite{langer-castelnuovo} do not work.  This
different proof is required by the nature of twisted sheaves and is
very similar to the proof included in \cite{twisted-moduli}.  One
amusing aspect of the theory of twisted sheaves that makes it hard to
simply carry over classical results is the fact that there is no
``trivial sheaf''; in fact, one does not even know \emph{a priori\/}
that there is a single semistable sheaf!  This is proven below in
Section \ref{sec:non-emptiness}.

The third difference lies in our independence from the methods of
Geometric Invariant Theory (GIT).  As we show, the standard O'Grady outline
really leads to a description of the stack of $\mu$-semistable sheaves
and its open substack of $\mu$-stable sheaves.  No mention of Gieseker
stability or the subtleties of GIT are necessary.  We strongly believe
that this approach clarifies the subject.  

Finally, let us note that while it may be tempting to view the present
work as a mere exercise in transferring classical results to the
orbifold context, this is deeply misguided.  The uninitiated reader is
referred to \cite{period-index-paper} for examples of applications of
the ideas discussed here and in \cite{twisted-moduli} to division
algebras over surfaces.  Adding an orbifold structure to the
underlying surface, giving rise to the spaces described in the present
work, seems to have some relation to the cyclicity problem for
division algebras of prime degree. While we have chosen to omit such
arithmetic considerations from this paper, they remain a strong
motivation for studying these moduli spaces.

\section*{Acknowledgments}
\label{sec:acknowledgments}

During the course of this work the author had helpful conversations
with numerous people, including Jean-Louis Colliot-Th\'el\`ene, Brian
Conrad, Johan de Jong, Joseph Lipman, Davesh Maulik, and Amnon Neeman.
During the preparation of the manuscript, the author was partially
supported by an NSF Postdoctoral Fellowship and NSF grant DMS-0758391.

\section{Notation}
\label{sec:notation}

We will say that a gerbe is \emph{non-trivial\/} if its structure
group is not the singleton group, and \emph{non-neutral\/} if it does
not have a section.  (Similar conventions hold with the ``non''s
removed!)  Thus, a residual gerbe on a stack is non-trivial when the
corresponding ``point'' is a stacky point, but one can certainly have
neutral non-trivial residual gerbes (e.g., if the field of moduli is
algebraically closed!).  

Let $X$ be a smooth proper geometrically connected orbisurface over a
field $k$.  By this we mean that $X$ has trivial inertia in
codimension 1, so that the stacky locus of $X$ consists of isolated
non-trivial residual gerbes.  We will write $X^\sp$ for the open
substack over which the inertia stack is trivial (the ``spatial
locus''); this is the largest open substack which is isomorphic to an
algebraic space.  Write $\widebar X$ for the coarse moduli space of
$X$ and $\sigma:X\to\widebar X$ for the map.

We will fix a $\m_N$-gerbe $\pi:\ms X\to X$ with $N$ invertible in $k$ for
the entire paper.  We will write $\D(\ms X)$ to denote the derived
category of quasi-coherent $\ms X$-twisted sheaves; we will use a
similar notation for any $\m_N$-gerbe.

By definition, given a morphism of finite type $X\to T$ of
stacks, a ``flat family of coherent sheaves over $T$'' is a
$T$-flat quasi-coherent $\ms O_X$-module of finite presentation.  (The
fibers become coherent, but coherent is not the correct property for
the family, as it is not preserved under base change on $T$.)

Given a smooth surface $Y$ over a field $k$, we will let $L_Y$ be the
invariant defined by Langer in section 2.1 of
\cite{langer-castelnuovo}.
\section{Preliminaries on orbisurfaces and regular sheaves}
\label{sec:prel-orbis}

\subsection{Intersection theory}
\label{sec:intersection-theory}

\begin{thm}\label{T:serre-duality}
  There exists a trace map $\H^2(X,\omega_X)\to k$ such that for all
  perfect complexes of quasi-coherent sheaves $\ms F$ and $\ms G$ on $X$ and any integer $i$, the pairing
$$\hom_{\D(X)}(\ms F,\ms G)\times\hom_{\D(X)}(\ms G,\ms F\tensor\omega_X[2])\to k$$ 
induced by the cup product and the trace maps is perfect.
\end{thm}
\begin{proof} For the proof, we refer the reader to Proposition 1.9,
  Proposition 1.14, and Theorem 1.32 of \cite{nironi}.
\end{proof}

\begin{defn}
  An invertible sheaf $\ms L$ on $X$ is \emph{ample\/} if the
  non-vanishing loci of sections of powers of $\ms L$ generate the
  topology on the underlying topological space $|X|$.  If $X$ has an
  ample invertible sheaf, we will say that $X$ is \emph{pseudo-projective\/}.
\end{defn}
Note that this definition is equivalent to the statement that some
power of $\ms L$ is the pullback of an ample invertible sheaf from
$\widebar X$.  There are more refined notions of ample sheaves
(related to embeddings in weighted projective stacks), but they will
have no place here (and seem to be overly specialized).

We use Vistoli's intersection theory on $X$ in what follows \cite{vistoli}.  In
particular, we have a theory of Chern classes and Todd classes, and a
degree map $\deg:A^2(X)\to\Q$.

\begin{thm}\label{T:pre-hodge}
  Suppose $H$ is an ample Cartier divisor on $X$.  If $D$ is a divisor
  such that $D^2>0$ then $D\cdot H\neq 0$.
\end{thm}
\begin{proof}
  The proof is standard and familiar from the theory of surfaces.
  First, note that To\"en's Riemann-Roch formula \cite{toen} says that
  $\chi(nD)$ is quadratic in $n$ with leading coefficient $D^2$.
  Thus, using Theorem \ref{T:serre-duality} and letting $K$ denote a
  canonical divisor on $X$, we have that either
  $h^0(nD)$ or $h^0(K-nD)$ goes to infinity with $n$, and similarly
  for either $h^0(-nD)$ or $h^0(K+nD)$.  

  \begin{claim}
    We cannot have both $h^0(K-nD)$ and $h^0(K+nD)$ going to infinity
    with $n$.
  \end{claim}
Indeed, given $\sigma\in\Gamma(X,K-nD)$, the map
$\tau\mapsto\tau\tensor\sigma:\Gamma(X,K+nD)\to\Gamma(X,2K)$ is
injective.  But then $h^0(K+nD)$ is bounded for all $n$ by $h^0(2K)$,
so it cannot grow with $n$.

We conclude that either $h^0(nD)$ or $h^0(-nD)$ is non-zero for some
$n>0$.  Since $H$ is ample, we have that $H\cdot(nD)>0$ or
$H\cdot(-nD)>0$.  In particular, $H\cdot D\neq 0$.
\end{proof}

\begin{cor}
  If $X$ is pseudo-projective then the intersection form on
  $NS(X)\tensor\R$ has signature $(1,r-1)$, the first
  basis vector coming from an ample divisor on $X$.
\end{cor}
\begin{proof}
  Diagonalizing the form over $\R$, the result follows from Theorem
  \ref{T:pre-hodge}.
\end{proof}

\subsection{Uniformizations of $\m_N$-gerbes on smooth
  pseudo-projective orbisurfaces}
\label{sec:uniformization}

We work with $\m_N$-gerbes throughout this paper.  For basic
properties of gerbes in the context of twisted sheaves, the reader is
referred to \cite{period-index-paper} and \cite{twisted-moduli}.  For
a proof that a gerbe on an orbisurface (or, more generally, a stack on
a stack) is an (algebraic) stack, the reader is referred to section
2.4 of \cite{paiitbgoaas}.

Let $k$ be an infinite field and $X$ a proper smooth geometrically connected
pseudo-projective orbisurface and $\ms X\to X$ a $\m_N$-gerbe, with
$N$ invertible in $k$.  
By standard results \cite{vistoli-kresch-etc}, \cite{vistoli-kresch}, \cite{dejong-gabber}, we know that $\ms X$ is a quotient stack. 
We recall the following result of Kresch and Vistoli \cite{vistoli-kresch}.

\begin{prop}
  There is a finite flat generically separable morphism $f:Y\to \ms X$
  with $Y$ a smooth geometrically connected projective surface over
  $k$.
\end{prop}
\begin{proof}[Idea of proof.]
  We indicate the idea of the proof so that we can point out why the
  cover is generically separable (a statement which was left out of
  \cite{vistoli-kresch}).  
The proof proceeds by using the fact that $\ms X$ is a quotient stack to
  note that there is a product $P\to\ms X$ of projective bundles which
  has a dense open representable substack whose complement has
  arbitrarily high codimension.  The coarse space $U$ of $P$ is then
  shown to be projective, and the cover $Y$ is found by taking general
  hyerplane sections of $U$.  Since a general such section will be
  generically separable over $\widebar X$, we will have that $f$ is
  generically separable.
\end{proof}

We can use the uniformization $f$ to give a quick proof of
algebraicity of certain stacks of sheaves.  For proofs in more general
contexts, the reader is referred to \cite{twisted-moduli} and \cite{MR2233719}.
This more ad hoc approach suffices for the task at hand.

It is elementary that the fibered category $Sh(\ms X)$ of flat
families of coherent $\ms O_{\ms X}$-modules is a stack over $k$.
Because it is locally free, the uniformization $f:Y\to\ms X$ induces a
pullback morphism $f^{\ast}:Sh(\ms X)\to Sh(Y)$.

\begin{prop}\label{L:pullback-yummy-uniformization}
  With the above notation, the morphism $f^{\ast}$ is representable by
  separated schemes of finite presentation.
\end{prop}
\begin{proof}
  Let $Y^{(2)}=Y\times_{\ms X}Y$ and $Y^{(3)}=Y\times_{\ms
    X}Y\times_{\ms X}Y$.  The two projections $p,q:Y^{(2)}\to Y$ are
  finite and flat, as are the three projections $Y^{(3)}\to Y^{(2)}$
  and (therefore) the three projections $Y^{(3)}\to Y$.
  (It follows that $Y^{(2)}$ and $Y^{(3)}$ are projective, but there
  is no reason to believe that they are smooth or geometrically connected.)

  Let $\phi:T\to Sh(Y)$ be a $1$-morphism, corresponding to a flat
  family of coherent sheaves $\ms F$ on $Y\times T$.  Define a functor
  $\ms D$ on the category of $T$-schemes as follows: given a
  $T$-scheme $S\to T$, the set $\ms D(S)$ consists of isomorphisms
  $\psi:p_S^\ast\ms F_S\simto q_S^{\ast}\ms F_S$ which satisfy the
  usual cocycle condition on $Y^{(3)}$.  It is clear that $\ms D$ is a
  locally closed subfunctor of the functor $H=\uhom(p^\ast\ms F,q^\ast\ms F)$
  parametrizing homomorphisms between the two pullbacks of $\ms F$ to
  $Y^{(2)}$.  Since $q^\ast\ms F$ is flat over $T$, it is well-known
  that $H$ is a $T$-scheme of finite presentation.  (That
  $\ms D$ is locally of finite presentation is immediate.)
\end{proof}
\begin{cor}
  The stack of coherent sheaves on $\ms X$ is an Artin stack locally
  of finite type over $k$.
\end{cor}
\begin{proof}
  By Proposition \ref{L:pullback-yummy-uniformization}, this follows
  from the corresponding fact for $Sh(Y)$, which is Th\'eor\`eme
  4.6.2.1 of \cite{l-mb}.
\end{proof}
\begin{cor}
  The stack of flat families of coherent $\ms X$-twisted sheaves is an
  Artin stack locally of finite type over $k$.
\end{cor}
\begin{proof}
  This follows immediately from the fact that this is an open (and
  closed) substack of the stack of coherent sheaves on $\ms X$.
\end{proof}

We can also use $f$ to define numerical conditions which yield
quasi-compact stacks of sheaves.  The polarization $H$ on $\widebar X$
induces a polarization $H|_Y$ on $Y$ which we will always use to define the
projective structure on $Y$.

\begin{defn}
  The \emph{$f$-Hilbert polynomial\/} of a coherent sheaf $\ms F$ on
  $\ms X$ is the Hilbert polynomial of $f^{\ast}\ms F$, i.e., the
  numerical polynomial $P^f_{\ms F}(m)=\chi(Y,\ms F_Y(mH_Y))$.
\end{defn}
\begin{lem}\label{L:rel-hilb-fam}
  The $f$-Hilbert polynomial is constant in a flat family of coherent
  sheaves on $\ms X$ parametrized by a connected base $T$.
\end{lem}
\begin{proof}
  This follows immediately from the corresponding
  well-known fact on $Y$.
\end{proof}

\begin{cor}\label{sec:unif-m_n-gerb}
  Given a coherent $\ms X$-twisted sheaf $\ms E$ and a numerical polynomial
  $P$, the scheme  $\Quot(\ms E,P)$ parametrizing flat families of
  coherent quotients of $\ms E$ with $f$-Hilbert polynomial is proper.
\end{cor}
\begin{proof}
  Quasi-compactness (and separatedness) follows immediately from Proposition
  \ref{L:pullback-yummy-uniformization} and the corresponding fact for
  the $\Quot$-scheme relative to $Y$.  Properness follows by the usual
  proof of the valuative criterion (see e.g.\ the proof of Theorem
  2.2.4 of \cite{h-l}).
\end{proof}

Further uses of the uniformization $f$ in studying semistable sheaves
will be given in section \ref{sec:moduli-slope-stable}.  

\subsection{$A_{N-1}$-orbisurfaces and the To\"en-Grothendieck-Hirzebruch-Riemann-Roch formula}
\label{sec:TGHRR}

In this section we discuss the orbisurfaces of interest to us and
compute the correction terms to the classical Hirzebruch-Riemann-Roch
formula for computing Euler characteristics of coherent sheaves.

A succint description of To\"en's Riemann-Roch formula can be found in
Appendix A of \cite{tseng} or \cite{toen}.  We refer the reader there
or to \cite{toen} for details.  The basic insight which To\"en
had is that the correct cohomology in which to work to produce a
Riemann-Roch formula on an orbifold is that of the inertia stack.
This yields correction terms to the standard Riemann-Roch formula
arising from the additional components of the inertia stack.  Our
short-term goal is a calculation of these correction terms.

\subsubsection{$A_{N-1}$-orbisurfaces}
\label{sec:a_n-1-orbisurfaces}

We recall a definition from \cite{paiitbgoaas}.

\begin{defn}
  An orbisurface $Z$ with coarse moduli space $\widebar Z$ is a
  \emph{Zariski $A_{N-1}$-orbisurface\/} if 
  \begin{enumerate}
  \item $N$ is invertible on $Z$;
  \item it has isolated non-trivial closed residual gerbes
    $\xi_1,\ldots,\xi_s$ with images $p_1,\ldots,p_s$ in $\widebar Z$;
  \item for each $i=1,\ldots,s$, the fiber product $Z\times_{\widebar
      Z}\spec\ms O_{\widebar Z,p_i}$ is isomorphic to the stack
    quotient $[\spec R/\m_N]$, where $R$ is a regular local ring of
    dimension $2$ with an action of $\m_N$ whose induced representation on the
    tangent space decomposes as $\chi\oplus\chi^{-1}$, where
    $\chi:\m_N\to\G_m$ is the natural character.
  \end{enumerate}
\end{defn}
Since we will only work with Zariski $A_{N-1}$-orbisurfaces in this
paper, we will refer to these simply as $A_{N-1}$-orbisurfaces.

If $\widebar Z$ is a normal projective surface with isolated
$A_{N-1}$-singularities (i.e., with local model
$k\[x,y,z\]/(z^n-xy)$), then it is the coarse moduli space of a smooth
pseudo-projective $\m_N$-orbisurface, unique up to unique isomorphism
over $\widebar Z$.  In our case the surface $\widebar Z$ will arise by
forming the cyclic cover of degree $N$ of a smooth projective surface
branched over a snc divisor.  In this case the singularities have the
form described above locally in the Zariski topology.

\subsubsection{The correction terms}
\label{sec:correction-terms}
In this section we compute the correction to the na\"ive Riemann-Roch
formula coming from the non-trivial residual gerbes of $X$.
\begin{lem}
  Let $\zeta$ be a primitive $N$th root of unity and $j<N$ a
  non-negative 
  integer.  We have 
$$\sum_{i=1}^{N-1}\frac{\zeta^{ij}}{2-\zeta^i-\zeta^{-i}}=\frac{j(j-N)}{2}+\frac{N^2-1}{12}.$$
\end{lem}
\begin{proof} 
Write $S=\mu_N(\widebar k)\setminus\{1\}$ and $1-S=\{y | 1-y\in S\}$. The sum to be
evaluated is
$$\sum_{x\in S}\frac{x^j}{(1-x)+(1-x^{-1})}.$$
Letting $y=1-x$, we can rewrite the sum as 
$$-\sum_{y\in 1-S}\frac{(1-y)^{j+1}}{y^2}.$$
Moreover, we know that the set over which the sum is taken is
precisely the set of roots of the polynomial
\begin{equation}\label{Eq:charpoly}
\sum_{t=0}^{N-1}\binom{N}{t+1}z^t=0.
\end{equation}  
Expanding the sum to be
evaluated yields the double sum 
$$-\sum_{y\in 1-S}\sum_{s=0}(-1)^s\binom{j+1}{s}y^{s-2}.$$

It thus suffices to evaluate the sums $\tau(r):=\sum_{y\in S-1} y^r$
for integers $r$ between $-2$ and $N-1$.  First suppose $r>0$; the sum
becomes 
$$\tau(r)=\sum_{x\in
  S}\sum_{q=0}^r(-1)^{r-q}\binom{r}{q}x^q=N$$ by an elementary
computation.  It is also clear that $\tau(0)=N-1$.  To compute $\tau(-1)$
and $\tau(-2)$, we argue as follows: letting $y_1,\ldots,y_{N-1}$ be
the elements of $S-1$, we see that
$$\tau(-1)=\sigma_{N-2}(y_1,\ldots,y_{N-1})/\sigma_{N-1}(y_1,\ldots,y_{N-1}),$$
where $\sigma_i$ denotes the $i$th elementary symmetric function in
$N-1$ variables.  Similarly,
$$\tau(-2)=\tau(-1)^2-2\sigma_{N-3}(y_1,\ldots,y_{N-1})/\sigma_{N-1}(y_1,\ldots,y_{N-1}).$$  
Using
\eqref{Eq:charpoly}, we see that
$\sigma_i(y_1,\ldots,y_{N-1})=\binom{N}{N-i}$.  Putting this together,
we find
$$\tau(-1)=\frac{N-1}{2}$$
and
$$\tau(-2)=\frac{(N-1)(5-N)}{12}.$$
Basic algebraic manipulations now yield
$$-\sum_{s=0}^{j+1}(-1)^s\binom{j+1}{s}\tau(s-2)=\frac{j(j-N)}{2}+\frac{N^2-1}{12},$$
as desired.
\end{proof}

Define a function $f:\Z/N\Z\to\Q$ by the formula 
$$f(\widebar x)=\frac{x(x-N)}{2}+\frac{N^2-1}{12},$$
where $0\leq x<N$ is one of the canonical representatives for
$\Z/N\Z$.  

Now suppose $F$ is a perfect complex on $X$.  Using the
To\"en-Riemann-Roch formula, we can write  
$$\chi(F)=\deg(\chern(F)\cdot\td_X)+\sum_{i=1}^n\delta_i(F),$$
where $\delta_i(F)$ is a correction term (a priori lying in
$\Q(\m_N)$, but actually lying in $\Q$) coming from contributions at
the residual gerbe $\xi_i$.  Write
$[\LLL\iota_i^{\ast}F]=\sum_{j=0}^{N-1}e^{(i)}_j\chi^j$ for the class in
$K$-theory of the derived fiber of $F$ over $\xi_i$.

\begin{lem}
  We have $\delta_i(F)=\frac{1}{N}\sum_{j=0}^{N-1}f(j)e^{(i)}_j$.
\end{lem}
\begin{proof}
  It suffices to prove this when $k$ is algebraically closed (as the
  Euler characteristic is invariant under field extension, as is the
  intersection theory).  In this case, the inertia stack $\ms I(X)\to
  X$ has the form
$$\ms I(X)=X\sqcup \sqcup_{i=1}^n \sqcup_{\ell=1}^{N-1}\B\m_N,$$
with one copy of $\sqcup_{\ell=1}^N\B\m_N$ mapping to each $\xi_i$ in
the natural way (via the identification provided by $\iota_i$).

Consider the copy of $\sqcup_{\ell=1}^{N-1}\B\m_N$ in the inertia
stack lying over $\xi_i$.  On the $\ell$th component, a section acts
on $\LLL\iota_i^{\ast}F$ via the $\ell$th power of the natural action.
Thus, decomposing $\sum_{j=1}^{N-1}e^{(i)}_j\chi^j$ into weighted
eigenbundles on the $\ell$th component yields
$\sum_{j=0}^{N-1}\zeta^{j\ell}e^{(i)}_j\chi^{j\ell}$, so that the
weighted Chern character $\widetilde{\chern}(F|_{\B\m_N})$ equals
$\sum_{j=0}^{N-1}\zeta^{j\ell}e^{(i)_j}$.

On the other hand, the Todd class which intervenes in the
To\"en-Riemann-Roch formula has the following form.  Let $E$ be the
locally free sheaf $\ms T_X|_{\B\m_N}$ on $\B\m_N$, viewed as the
$\ell$th copy over $\xi_i$.  We know that
$E=\chi^{\ell}\oplus\chi^{-\ell}$, since $X$ is a 
$A_{N-1}$-orbisurface.  Using the standard formula (see
e.g. A.0.5 of \cite{tseng}), we find that the contribution to the To\"en-Todd class
coming from this component is then  
$$\widetilde{\td}(E)=\frac{1}{2-\zeta^{\ell}-\zeta^{-\ell}}.$$
Thus, taking into account the fact that $\B\m_n\to\spec\kappa$ has
degree $1/N$, we find that the total correction term for the
components of $\ms I(X)$ lying over $\xi_i$ is
$$\frac{1}{N}\sum_{\ell=1}^{N-1}\sum_{j=0}^{\N-1}\frac{\zeta^{j\ell}e^{(i)}_j}{2-\zeta^\ell-\zeta^{-\ell}}=\frac{1}{N}\sum_{j=0}^{N-1}e^{(i)}_j\sum_{\ell=1}^{N-1}\frac{\zeta^{j\ell}}{2-\zeta^{\ell}-\zeta^{-\ell}},$$
which proves the lemma.
\end{proof}
The formula provides an immediate corollary.
\begin{cor}\label{C:euler-structure-sheaf}
  The Euler characteristic of $\ms O_X$ satisfies
$$\chi(X,\ms O_X)=\deg\td_X+n\frac{N^2-1}{12N}.$$
\end{cor}
\subsection{Twisted sheaves on $A_{N-1}$-orbisurfaces}
\label{sec:sheav-deform-theory}

In this section we fix an $A_{N-1}$-orbisurface $X$ with coarse moduli
space $\widebar X$ and a $\m_N$-gerbe $\ms X\to X$. We assume that the
base field $k$ is algebraically closed and denote the non-trivial
residual gerbes of $X$ by $\xi_1,\ldots,\xi_n$, each of which is
isomorphic to $\B\m_N$ via a map $\iota_i:\B\m_N\inj X$.  Given
$i=1,\ldots,n$, we write $\ms X_i=\ms X\times_X\xi_i$.  Fix
$(-1)$-fold invertible $\ms X_i$-twisted sheaves $\ms L_i$.  Twisting
by $\ms L_i$ and pulling back by $\iota_i$ defines an equivalence of
categories between $\D(\ms X_i)$ and the category of graded
representations of $\m_N$.

\begin{notation}\label{N:otation}
Given a perfect complex $F$ of $\ms
X$-twisted sheaves, we thus get a class in $K(\m_N)$  by taking
the class in $K$-theory associated to the perfect complex $\ms
L_i\tensor\LLL\iota_i^{\ast}F$ of $\m_N$-representations.  We will
always write $\overline F_i$ for this class.
\end{notation}

We establish some of the basic properties of torsion free $\ms
X$-twisted sheaves which are relevant to the problems at hand; in
particular, we will focus on the role played by the non-trivial
residual gerbes of $X$.  The reader unfamiliar with twisted sheaves is
referred to \cite{twisted-moduli} and \cite{period-index-paper} (and
the references therein); we will not redevelop the theory from scratch in
the present paper.

While we write in the language of twisted sheaves, when $\ms X$ is the
trivial $\m_N$-gerbe, our results also apply to ordinary sheaves (and
perfect complexes thereof) on $X$.  Removing the twisting class is
straightforward in this case, but we will occasionally point out
explicitly how a certain definition or property would look in the
untwisted case.

Before proceeding, we record an orphan lemma, which is a twisted form
of Theorem \ref{T:serre-duality}.

\begin{lem}\label{sec:twisted-sheaves-a_n}
  Let $\ms E$ and $\ms H$ be perfect complexes of quasi-coherent $\ms
  X$-twisted sheaves.  The cup product and trace maps induce a perfect
  pairing
$$\hom(\ms E,\ms H)\times\hom(\ms H,\ms E\tensor\omega_X[2])\to k.$$
\end{lem}
\begin{proof}
  This follows from Theorem \ref{T:serre-duality} with $\ms F=\ms O_X$
  and $\ms G=\rshom(\ms F,\ms G)$, using the standard properties of
  perfect complexes.
\end{proof}

We start by noting a convention which will prove helpful.
\begin{defn}\label{sec:sheav-deform-theory-2}
  Given a perfect complex $F$ of $\ms X$-twisted sheaves, the
  \emph{normalized Chern classes\/}, denoted $\widetilde c_i(F)$,
  are $Nc_i(F)$, where $c_i(F)$ denotes the usual $i$th Chern class in
  $A^i(\ms X)$.  
\end{defn}

The purpose of normalizing the Chern classes is simply to correct
the fact that $\pi_{\ast}\pi^{\ast}$ acts as multiplication by $1/N$
on $c_i$.  Thus, if $G$ is a perfect complex on $X$, we have
that $\pi_{\ast}\widetilde c_i(\pi^{\ast}G)=c_i(G)$.

\begin{defn}\label{D:sing-locus}
  Given an $\ms X$-twisted sheaf $\ms F$, write $\sing(\ms F)$ for the reduced
  closed substack structure on the complement of the locus over which
  $\ms F$ is locally free.
\end{defn}
Since $\ms X$ and $X$ have the same topology, we can consider
$\sing(\ms F)$ as a subset of $X$.  (This serves primarily to simplify
notation and we will freely avail ourselves of this fact in the future.)

Since $k$ is algebraically closed, we can define the length of a
$0$-dimensional twisted sheaf as follows.
\begin{defn}\label{D:length}
  Given an $\ms X$-twisted sheaf $\ms Q$ of dimension $0$, the
  \emph{length\/} of $\ms Q$ is the maximal length of a filtration
  $\ms Q=\ms Q_{\ell}\supset\ms Q_{\ell-1}\supset\cdots\supset\ms
  Q_0=0$ with non-trivial subquotients.
\end{defn}

Using the convention of Definition \ref{sec:sheav-deform-theory-2},
one can see that $\ell(\ms Q)=-\widetilde c_2(\ms Q)$.

\begin{defn}\label{D:tot-reg}
  A perfect complex $\ms F\in\D(\ms X)$ is \emph{totally regular\/} if
  $\overline{\ms F}_i$ is in the ideal $\Z\rho\subset K(\m_N)$ for
  each $i=1,\ldots,n$.
\end{defn}

A weaker notion which will often be sufficient for certain statements
(but which will be relatively unimportant for this work) is the
following.

\begin{defn}\label{D:tot-pos}
  A perfect complex $\ms F\in\D(\ms X)$ is \emph{totally positive\/} if
  $\overline{\ms F}_i$ is the class in $K(\m_N)$ associated to a
  (non-virtual) representation of $\m_N$ for each $i=1,\ldots,n$.
\end{defn}

The reader can readily check that these two properties are independent
of the choices of the invertible sheaves $\ms L_i$ on $\ms X_i$ made
at the beginning of the section.

Since any coherent sheaf on $X$ is perfect as an object of the derived category
($X$ being regular), we can ascribe the same properties to sheaves.
We will do this freely.

\begin{lem}\label{L:tot-reg-props}
  If $\ms F$ is a totally regular perfect complex of $\ms X$-twisted
  sheaves and $\ms G$ is an arbitrary perfect complex then
  \begin{enumerate}
  \item $\ms F\ltensor\ms G$ is totally regular;
  \item $\rhom(\ms F,\ms G)$ is totally regular.
  \end{enumerate}
\end{lem}
\begin{proof}
  This immediately reduces to the analogous statement for objects in
  the derived category of $\m_N$-representations, where this follows
  from the fact that $\Z\rho$ is an ideal of $K(\m_N)$.  (We can use
  the tensor powers of the $\ms L_i$ on the various powers of $\ms X$
  to get complexes of $\m_N$-representations out of $\rshom(\ms F,\ms
  G)$ or $\ms F\ltensor\ms G$; since total regularity is independent
  of this local trivialization, we lose nothing by leaving such minor
  details unremarked upon.)
\end{proof}

Thus, if $\ms F$ is a totally regular perfect complex of $\ms
X$-twisted sheaves, we can identify (via derived pushforward) the
complex $\rshom(\ms F,\ms F)$ with a totally perfect complex on $X$.
We will implicitly do this below.

One thing which the reader should note is that the reflexive hull of a
totally regular torsion free $\ms X$-twisted sheaf on $X$ is not
necessarily totally regular.  It is a simple matter to make local
examples around any singular point, an exercise which we leave to the
reader.

\begin{prop}\label{P:tot-reg-rr}
  If $\ms F$ is a totally regular perfect complex on $X$ then 
$\chi(X,\ms F)=\deg(\chern(\ms F)\cdot\td_X).$
\end{prop}
\begin{proof}
This follows immediately from the calculations of section
\ref{sec:correction-terms}, together with the fact that
$\sum_{j=0}^{N-1}f(j)=0$.
\end{proof}

In what follows we use the fact that for any field $L$ there is a
canonical ring isomorphism $K(\m_{N,L})=\Z[x]/(x^N-1)$ (which is functorial
in $L$).  We will write $K(\m_N)$ for this constant ring.

\begin{prop}\label{P:repn-invariant}
  Let $T$ be a connected scheme and $\ms F$ a perfect complex on
  $T\times\B\m_N$.  There exists class $c\in K(\m_N)$ such that for
  every geometric point 
  $t:\spec\kappa\to T$, the class of the preimage $\LLL
  t^{\ast}\ms F$ in $K(\m_N)$ is $c$.
\end{prop}
\begin{proof}
  We may assume that $T$ is affine, so that there is a global
  resolution of $\ms F$ by locally free sheaves.  This reduces us to
  the case in which $\ms F$ is itself locally free.  The result now
  follows either from the theory of Hopf comodules for the Hopf
  algebra of $\m_N$, or from the upper semicontinuity of cohomology
  for flat coherent sheaves on proper morphisms of Artin stacks.
\end{proof}

The deformation theory of $\ms X$-twisted sheaves is governed by
deformations of sheaves of modules in the \'etale topos of $\ms X$ and
is thus described by the general theory of Illusie.  In
particular, we have the following proposition.  Given a perfect
complex $\ms F$ and an ideal $I$, there is a trace map $\ext^i(\ms
F,\ms F\ltensor I)\to\H^i(\ms X,I)$ arising from the trace map
$\rshom(\ms F,\ms F)\to\ms O$ and the identification $\ext^i(\ms F,\ms
F\tensor I)=\ext^i(\ms O,\rshom(\ms F,\ms F)\ltensor I)$.  We let
$\ext^i(\ms F,\ms F\ltensor I)_0$ denote the kernel of the trace map.
As usual, when the rank of $\ms F$ is invertible in $k$, Serre duality
induces an isomorphism $\ext^2(\ms F,\ms F)_0\simto\hom(\ms F,\ms
F\tensor\omega_X)_0^{\vee}$.

\begin{prop}\label{sec:sheav-deform-theory-1}
Let $A\to A_0$ be a surjection of $k$-algebras with kernel $I$ of square
$0$.  Let $\ms F$ be a flat family of coherent $\ms X$-twisted sheaves
parametrized by $A$ and let $\ms L_0$ be the determinant of $\ms F$.
Fix a flat family $\ms L$ of invertible sheaves on $\ms X$ parametrized by $A$
along with an isomorphism $\ms L\tensor_A A_0\simto\ms L_0$.
\begin{enumerate}
\item There is an element $\mf o$ of $\ext^2(\ms F,\ms F\ltensor I)_0$
  such that $\mf o=0$ if and only if there is a flat family of
  coherent $\ms X$-twisted sheaves $\F$ parametrized by $A$ and
  isomorphisms $\F\tensor_A A_0\simto\ms F$ and $\det\F\simto\ms L$
  compatible with the given isomorphism $\ms L\tensor_A A_0\simto\ms
  L_0$ and the identification of $\ms L_0$ with the determinant of
  $\ms F$.
\item The set of isomorphism classes of such extensions is a
  pseudo-torsor under $\ext^1(\ms F,\ms F\ltensor I)_0$.
\item The set of infinitesimal automorphisms of one such extension is trivial.
\end{enumerate}
\end{prop}

\begin{defn}\label{D:good}
  A coherent $\ms X$-twisted sheaf $\ms F$ is \emph{unobstructed\/} if
  $\ext^2(\ms F,\ms F)_0=0$.
\end{defn}

\begin{defn}
  Let $T$ be a $k$-scheme.  A \emph{$T$-flat family of torsion free
    coherent sheaves on $X$\/} is a quasi-coherent sheaf $\F$ on
  $X\times T$ which is locally of finite presentation, flat over $T$,
  and such that for every geometric point $t\to T$, the fiber $\ms
  F_t$ is a torsion free coherent sheaf on $X\tensor_k\kappa(t)$.
\end{defn}

\begin{defn}
  Given a torsion free sheaf $\ms F$ on $\ms X$,
  \begin{enumerate}
  \item the \emph{soft singular locus\/} of $\ms F$ is $\sing^\soft(\ms
    F):=\sing(\ms F)\cap X^\sp$;
  \item the \emph{soft hull\/} of $\ms F$ is the sheaf $\ms F^\soft$ with
    inclusion $\ms F\inj\ms F^\soft$ such that $(\ms F\inj\ms
    F^\soft)|_{X^\sp}=\ms F|_{X^\sp}\inj\ms F^{\vee\vee}_{X^\sp}$ and $\ms
    F\inj\ms F^\soft$ is isomorphic to the identity map in a neighborhood
    of $\sing(\ms F)\setminus\sing^\soft(\ms F)$;
  \item the \emph{soft colength\/} of $\ms F$ is $\ell^\soft(\ms
    F):=\ell(\ms F^\soft/\ms F)$.
  \end{enumerate}
\end{defn}

\begin{defn}
  A torsion free sheaf $\ms F$ on $X$ is \emph{softly reflexive\/} if
  the map $\ms F\to\ms F^\soft$ is an isomorphism.
\end{defn}

Given a field $L\subset k$ and a sheaf $\ms G$ on $X\tensor L$, the
length of $\ms G$ will always mean the length of $\ms G\tensor\widebar
L$ for any chosen algebraic closure of $L$.  Thus, when we speak of
the colength of the fiber of a family of torsion free sheaves, we mean
the colength of the geometric fiber.

\begin{prop}\label{P:soft-colength-constr}
  Suppose $\ms F$ is a $T$-flat family of torsion free 
  coherent sheaves on $\ms X\times T$. The function $t\mapsto\ell^\soft(\ms
  F_t)$ is constructible.
\end{prop}
\begin{proof}
  It suffices to prove the following: assuming that $T$ is integral,
  we have that $\ell^\soft(\ms F_\eta)=n$ if and only if there is a dense
  open subset $W\subset T$ such that $\ell^\soft(\ms F_w)=n$ for all $w\in
  W$.

The inclusion $\ms F_{\eta}\inj\ms F_{\eta}^\soft$ extends to a short
exact sequence
$$0\to\ms F_W\to\ms F_W'\to Q\to 0$$
over some open subset $W$ of $T$, where $\ms F_W'$ has softly
reflexive fibers.  Moreover, shrinking $W$ if
necessary, we may assume that $Q$ is flat and quasi-finite over $W$.
Finally, since the support of $Q_{\eta}$ is contained in $X^\sp_\eta$,
we see that further shrinking $W$ we may suppose that $\supp Q\subset
X^\sp\times W$.  By Zariski's Main Theorem, after shrinking $W$ one
more time we may assume that $Q$ is finite over $W$.  
It follows that each fiber of the exact sequence
computes the soft hull of the fibers of $\ms F$, and the degree of
$Q$ over $W$ computes the (constant) colength of the (geometric) fibers.
\end{proof}

\begin{defn}\label{D:soft-bdry}
  Given a $T$-flat family $\ms F$ of torsion free coherent $\ms
  X$-twisted sheaves parametrized by 
  $T$, the \emph{soft boundary\/} of the family is 
$${\softboundary}T:=\{t\in T | \ell^\soft(\ms F_t)>0\}.$$
\end{defn}

By Proposition \ref{P:soft-colength-constr}, we know that ${\softboundary}T$
is a locally constructible subset of $T$.  Given a quasi-compact open
subscheme $T'\subset T$, we can thus stratify ${\softboundary}T'$ by
locally closed subschemes.  In particular, there will be a largest
locally closed stratum, $({\softboundary}T')^\dom$.

\begin{lem}\label{L:bdry-codim}
  Suppose $T$ is locally Noetherian.  With the preceding notation, if
  ${\softboundary}T\neq\emptyset$, then for all quasi-compact open
  subschemes $T'\subset T$ we have $\codim(({\softboundary}T')^\dom,T')\leq
  r-1$.  Equivalently, if $t\in{\softboundary}T$ is a minimal point then
  $\dim\ms O_{T,t}\leq r-1$.
\end{lem}
\begin{proof}
  A standard argument (e.g., Lemma 9.2.1 of \cite{h-l}) shows that 
$$D:=\{(t,u) | \ms F_{(t,u)}\textrm{ is not locally free}\}$$
has codimension at most $r+1$ in $X^\sp\times T$.  Moreover, $D\to T$ is
quasi-finite and of finite type.  It follows that the minimal points
of the image all have codimension at most $r+1-2=r-1$, as desired.
\end{proof}

\subsection{Soft $\Quot$ schemes}
\label{sec:soft-quot-schemes-1}

Let $\ms F$ be a softly reflexive torsion free $\ms X$-twisted sheaf.

\begin{defn}
  The \emph{soft $\Quot$ space of length $\ell$ quotients\/} of $\ms
  F$, denoted $\Quot^{\soft}(\ms F,\ell)$, is the algebraic space
  parametrizing quotients $\ms F\to\ms Q$ with $\ms Q$ a coherent
  twisted sheaf of length $\ell$ with support contained in $X^\sp$.
\end{defn}

It is a standard application of Artin's theorem that this functor is
an algebraic space.  (The reader can see, e.g., \cite{olsson-starr} for a
proof.)  Since we constrain the support of the quotient to the spatial
locus of $X$, we cannot hope that $\Quot^{\soft}$ is proper.  However,
we do have the following irreducibility result.

\begin{prop}\label{sec:soft-quot-schemes-2}
  The space $\Quot^{\soft}(\ms F,\ell)$ is irreducible of dimension
  $\ell(\rk\ms F+1)$.
\end{prop}
\begin{proof}
  The proof is exactly analogous to the proof of Lemma 2.2.7.28 of
  \cite{twisted-moduli}, and reduces to the proof of the classical
  case by Ellingsrud-Lehn \cite{MR1714770}.
\end{proof}

\section{Stability and boundedness for twisted sheaves on $X$}
\label{sec:stab-bound}

\subsection{Slopes, $\mu$-stability, and the Segre invariant}
\label{sec:slop-segre-invar}

Let $H$ be a fixed ample divisor on the coarse space $\widebar X$.
(We remind the reader that we denote the natural map $X\to\widebar X$
by $\sigma$.)

\begin{defn}
  Let $F$ be a perfect complex on $\ms X$ with positive rank $r$.  The
  \emph{slope\/} of $F$ (\emph{with respect to $H$\/}) is 
$$\mu(F)=\frac{H\cdot\widetilde c_1(F)}{r}.$$
\end{defn}
It is easy to see that the slope of $F$ is a rational number with
denominator dividing $rN$.  In particular, any bounded set of
differences of slopes of perfect complexes on $\ms X$ has a minimal
element.

\begin{lem}
  If $F$ is a perfect complex on $\ms X\times T$ with $T$ a connected
  scheme then for all pairs of geometric points $s,t\to T$, the slopes 
$\mu(F_t)$ and $\mu(F_s)$ are equal.
\end{lem}
\begin{proof}
  The proof is straightforward and comes down to showing that
  intersections of invertible sheaves are constant in a connected
  family.  We refer the reader to Proposition 2.2.7.22 of
  \cite{twisted-moduli} for an analogous
  proof with details.
\end{proof}

\begin{lem}\label{L:slope-hilbert-comp}
  If $\ms E$ is a coherent sheaf on $X$ with Hilbert polynomial
  $$P(m)=\chi(X,\ms E(m))=\alpha_2(\ms E)\frac{m^2}{2}+\alpha_1(\ms
  E)m+\alpha_0(\ms E)$$ 
then 
$$\deg \ms E=\alpha_1(\ms E)-\alpha_1(\ms O_X)\rk\ms E.$$
In particular, the slope of $\sigma_{\ast}\ms E$ equals the slope of
$\ms E$.
\end{lem}
\begin{proof}
We remind the reader that the slope of a coherent sheaf on $\widebar
X$ is computed as in Definition 1.2.11 of \cite{h-l}, using the coefficients
of the Hilbert polynomial.  (The formula is the one indicated as the
conclusion of the lemma.)  

This follows immediately from the To\"en-Grothendieck-Riemann-Roch
formula.    
\end{proof}

\begin{defn}
  A torsion free sheaf $E$ on $\ms X$ is \emph{$\mu$-stable\/} if
  $\mu(F)<\mu(E)$ for all subsheaves $F\subset E$ such that $0<\rk
  F<\rk E$.  The sheaf $E$ is \emph{$\mu$-semistable\/} if
  $\mu(F)\leq\mu(E)$ for all non-zero subsheaves $F\subset E$.
\end{defn}
When $E$ is an $\ms X$-twisted sheaf, every subsheaf will also be $\ms
X$-twisted, so the notion of stability for a twisted sheaf depends
purely on the category of $\ms X$-twisted sheaves.  In addition, it is
standard that to verify semistability, it suffices to check that the
slope does not increase among the same set of subsheaves as is used
for checking stability.

One can measure the degree of stability of a sheaf using the following
invariant.

\begin{defn}
  Given a torsion free sheaf $E$ of rank at least $2$ on $\ms X$, the
  \emph{Segre invariant\/} of $E$, denoted $s(E)$, is the minimal
  value of $\mu(E_2)-\mu(E_1)$, where $E_1$ and $E_2$ are torsion free
  sheaves fitting into an exact sequence
$$0\to E_1\to E\to E_2\to 0.$$
\end{defn}

The existence and uniqueness of the maximal destabilizing subsheaf
(by a proof formally identical to that in \cite{h-l}) imply that the set of
such differences is in fact bounded and therefore there is a minimal
such value (as implicitly asserted in the definition).

The Segre invariant, originally defined and studied by Segre and
Nagata, was systematically studied by Langer in the context of moduli
of sheaves \cite{langer-castelnuovo}.  He observed that the Segre
invariant provides a satisfactory replacement for the notion of
$e$-stability used by O'Grady \cite{o'grady} when one wishes to study moduli in
a characteristic-agnostic manner.  Very positive Segre invariants
indicate a high degree of stability, while very negative Segre
invariants indicate a high degree of instability.

The following results will be useful to us; their proofs are formally
identical to those in Corollary 2.3 and Lemma 2.4 of \cite{MR2085175}, so we omit
them here.  (The fact that we work with sheaves on stacks and Langer
works with sheaves on varieties should not cause alarm -- only the
formal properties of the abelian categories of sheaves and the
intersection theory play a role in the proofs.)

\begin{lem}\label{L:slop-stab-segre-1}
  Let $E$ be a torsion free sheaf on $\ms X$ of rank at least $2$, and
  assume that $s(E)$ is realized by an exact sequence $0\to E_1\to
  E\to E_2\to 0$.
  \begin{enumerate}
  \item If $E$ is $\mu$-semistable then both $E_1$ and $E_2$ are
    $\mu$-semistable.
  \item If $E$ is $\mu$-stable then both $E_1$ and $E_2$ are $\mu$-stable.
  \end{enumerate}
\end{lem}

The following is proven more generally for higher-dimensional
varieties in Lemma 2.4 of \cite{MR2085175}, but we will only need it
for surfaces.

\begin{lem}\label{L:slop-stab-segre-2}
  Suppose $E$ is a torsion free sheaf on $\ms X$ of rank at least $2$.
  Let $D\subset X$ be a normal irreducible curve such that $E|_D$ is
  locally free.  If $F$ is an elementary transformation of $E$ along
  a quotient of $E|_D$ then $s(E)\leq s(F)+D\cdot H$.
\end{lem}

The following is referred to by Langer without proof.  We give a proof
here.

\begin{lem}\label{sec:slop-stab-segre}
  Let $s\in\Q$ and let $F$ be a flat family of torsion free sheaves on
  $\ms X\times T$.  There is a closed subset $T(s)\subset T$
  parametrizing the geometric points $t$ such that $s(F_t)\leq s$.
\end{lem}
\begin{proof}
  First, we argue that the set is locally constructible; we may assume
  that $T$ is Noetherian.  Theorem 2.3.2 of \cite{h-l} shows that the
  Harder-Narasimhan filtrations on the $F_t$ form a flat family over
  some dense open subscheme $U\subset T$.  It follows (by stratifying
  $T$) that there is a maximal value $m$ for the slope of a subsheaf of a
  (geometric) fiber of $F$, i.e., for all exact sequences
$$0\to E_1\to F_t\to E_2\to 0$$
we have that $\mu(E_1)\leq m$.  In this case, $\mu(E_2)-\mu(E_1)\geq
\mu(E_2)-m$, so that when $s(E_t)\leq s$ we have $\mu(E_2)\leq s+m$.
On the other hand, Grothendieck's boundedness theorem (Theorem 1.7.9
of \cite{h-l}) shows that the scheme $Q$ of flat families of quotients of $F$
whose fibers have slopes bounded above by $s+m$ is of finite type.  On
the scheme $Q$, we can look at the closed subscheme $Q'\subset Q$
parametrizing exact sequences as above with $\mu(E_2)-\mu(E_1)\leq
s$ (noting that the difference being computed is a locally constant
function on $Q$ with finite range).  By Chevalley's theorem, $T(s)$ is
constructible.

On the other hand, it is clear that $T(s)$ is closed under
specialization, as we can spread out (over a base dvr) any sequence
giving an upper bound for the Segre invariant (and the resulting
slopes are constant in the fibers over the base dvr).  Thus,
$T(s)$ is closed, as desired.
\end{proof}

There is one final basic lemma concerning endomorphisms of semistable
sheaves which we will need in the sequel.

\begin{lem}\label{L:simplest-bound}
  If $\ms F$ is a $\mu$-semistable $\ms X$-twisted sheaf of rank $r$
  then $\dim\hom(\ms F,\ms F)\leq r^2$.
\end{lem}
\begin{proof}
  Any endomorphism of $\ms F$ must preserve the socle (see Lemma
  1.5.5ff of \cite{h-l});
  moreover, the quotient $\ms F/\operatorname{Soc}(\ms F)$ is also
  semistable.  The result follows by induction from the polystable
  case, which itself follows immediately from the fact that stable
  sheaves are simple.
\end{proof}

\subsection{Discriminants and the Bogomolov inequality}
\label{sec:discr-bogom-ineq}

In this section we fix a uniformization $f:Y\to\ms X$ as in section
\ref{sec:uniformization}.  Let $d$ be the degree of $Y$
over $X$ (not over $\ms X$ -- so $[Y:\ms X]=Nd$), and define $L_X=L_Y/d$,
where $L_Y$ is the invariant defined by Langer in section 2.1 of \cite{langer-castelnuovo}.

\begin{defn}
Given a coherent $\ms X$-twisted sheaf $\ms E$, the discriminant of
$\ms E$, denoted $\Delta(\ms E)$, is defined by
$$\Delta(\ms E):=2\rk(\ms E)\widetilde c_2(\ms E)-(\rk\ms
E-1)\widetilde c_1^2(\ms E).$$
\end{defn}
It is elementary that $\Delta(\ms E)=c_2(\R\pi_{\ast}\R\shom(\ms E,\ms
E))$, and that the discriminant is invariant in a flat family of
$\ms X$-twisted sheaves over a
connected base.  

\begin{notation}
  Given $L\in\Pic(\ms X)$, we will write $\Tw_{\ms
    X}^{ss}(r,L,\Delta)$ for the stack of pairs $(\ms F,\psi)$, where
  $\ms F$ is a totally regular $\mu$-semistable $\ms X$-twisted sheaf of
  rank $r$ with discriminant $\Delta$ and $\psi:\det\ms F\simto L$ is
  an isomorphism.
\end{notation}

\begin{lem}\label{sec:discr-bogom-ineq-2}
  If $\ms E$ is totally regular of rank $r$ then
$$\chi(\ms E,\ms E):=\chi(\rshom(\ms E,\ms E))=-\Delta+r^2\left(\chi(\ms O_X)-n\frac{N^2-1}{12N}\right).$$
\end{lem}
\begin{proof}
  Since $\rshom(\ms E,\ms E)$ is totally regular with trivial
  determinant and rank $r^2$, Proposition \ref{P:tot-reg-rr} shows
  that
$$\chi(\ms E,\ms E)=-\Delta+r^2\deg\td_X.$$
The formula now follows from Corollary \ref{C:euler-structure-sheaf}.
\end{proof}

As we will discuss below in section \ref{sec:dimension-estimates}, the
dimension of $\Tw_{\ms X}^{ss}(r,L,\Delta)$ depends linearly upon
$\Delta$.

Given a torsion free coherent $\ms X$-twisted sheaf $F$, the
normalizations we have established show that $\mu(f^{\ast}F)=d\mu(F)$
and $\Delta(f^{\ast}F)=d\Delta(F)$.

\begin{prop}[Orbifold twisted Langer-Bogomolov
  inequality]\label{sec:discr-bogom-ineq-1}
  If $E$ is a semistable $\ms X$-twisted sheaf then
$$\Delta(E)\geq -\frac{L_X^2}{4H^2}\rk(E)^2(\rk(E)-1).$$
\end{prop}
\begin{proof}
  This follows immediately from the statement of Theorem 2.2 of
  \cite{langer-castelnuovo} applied to $f^\ast E$, using the fact that $f^\ast
  F$ is semistable (see Lemma \ref{L:i-cover-you} below).
\end{proof}

\subsection{Moduli of $\mu$-semistable twisted sheaves}
\label{sec:moduli-slope-stable}

We retain the choice of uniformization $f:Y\to\ms X$
of section \ref{sec:discr-bogom-ineq}.

In this section, we remark on the stack of semistable $\ms
X$-twisted sheaves.  Since we are not concerned with GIT quotients, we
can consider the stack of all $\mu$-semistable sheaves and do not need
to restrict our attention to an open sublocus of ``Gieseker semistable
sheaves.''  The reader interested in such a thing can see section
2.2.7.5 of \cite{twisted-moduli} for
one possible version of Gieseker stability in this context; the
methods used there seems likely to carry over to the present
situation.

The great advantage provided by considering only $\mu$-stability is
that it affords an especially efficient development of the theory in
the case at hand (i.e., a surface over a field).  The reason for this
is the following well-known lemma.

\begin{lem}\label{L:i-cover-you}
  A torsion free $\ms X$-twisted sheaf $\ms F$ is $\mu$-semistable if
  and only if $f^{\ast}\ms F$ is a $\mu$-semistable sheaf on $Y$.
\end{lem}
\begin{proof}
  The proof is identical to the proof of Lemma 3.2.2 of \cite{h-l}.  
\end{proof}

The reader should note that the same is not true for $\mu$-stability.

This covering is useful for proving a basic algebraicity and
boundedness statement.  More general versions of this statement are
true (cf.\ section 2.3.2 of \cite{twisted-moduli}), but this is all we
will need in this paper.

\begin{lem}\label{sec:moduli-mu-semistable-2}
  The stack of $\mu$-semistable $\ms X$-twisted sheaves with fixed
  $f$-Hilbert polynomial is of finite type over $k$.
\end{lem}
\begin{proof}
  This follows immediately from Lemma
  \ref{L:pullback-yummy-uniformization} and the corresponding fact for
  smooth projective surfaces (Theorem 4.2 of \cite{langer}).
\end{proof}

\begin{prop}\label{sec:moduli-mu-semistable-1}
  Given $L\in\Pic(\ms X)$, $r\geq 0$, and $\Delta\geq 0$, the stack
  $\Tw^{ss}_{\ms X}(r,L,\Delta)$ an Artin stack of finite type over
  $k$.
\end{prop}
\begin{proof}
Given $\ms F$, we know that $\Delta(f^{\ast}\ms F)=d\Delta(\ms F)$.
On the other hand, fixing the rank, determinant, and discriminant of a sheaf
on $Y$ fixes its Hilbert polynomial (by appeal to the Riemann-Roch
formula).  Thus, the sheaves parametrized by $\Tw^{ss}_{\ms
  X}(r,L,\Delta)$ have constant $f$-Hilbert polynomial.  The result
follows from Lemma \ref{sec:moduli-mu-semistable-2}.
\end{proof}
An alternative approach to proving the results of Proposition
\ref{sec:moduli-mu-semistable-1} in the relative case (where a finite
flat uniformization $Y\to\ms X$ may not exist) is given in section
2.3 of \cite{twisted-moduli}.  The algebraicity is easily proven
using Artin's theorem, while boundedness is a bit more complicated but
not very difficult.

\subsubsection{Pushing forward: openness of stability}
\label{sec:comp-stab-noti}

We now turn to the relationship between stable and semistable sheaves.
Again, in the case at hand there is a simple way to show that
$\mu$-(semi)stability is open in a flat family.  (For a proof in the
relative case which uses different techniques and which is easily
adaptable to the orbifold situation, the reader is again referred to
Corollary 2.3.2.12 of \cite{twisted-moduli}.)  We let
$\Tw_{\ms X}$ denote the stack of coherent $\ms X$-twisted sheaves.

Let $\ms V$ be a locally free $\ms X$-twisted sheaf of positive rank
(for exampe, the reflexive hull of any coherent $\ms X$-twisted sheaf
of positive rank, which will be locally free because $X$ is regular
and $2$-dimensional).  Let $\ms B=\pi_{\ast}\sigma_{\ast}\send(\ms
V)$; this is a coherent sheaf of algebras on $\widebar X$.  An $\ms
X$-twisted sheaf $\ms F$ gives rise to a right $\ms B$-module $M(\ms
F):=\pi_{\ast}\sigma_{\ast}\shom(\ms V,\ms F)$; this defines a map
from $\Tw_{\ms X}$ to the stack $\ms M_{\ms B}$ of $\ms B$-modules
(since $\pi_{\ast}$ takes flat families to flat families, $X$ being
tame).  We will write this morphism as $M:\Tw_{\ms X}\to\ms M_{\ms
  B}$.

Given a coherent sheaf $\ms C$ of $\ms O$-algebras (on either $X$ or
$\widebar X$), we will let $\Coh_{2,1}(\ms C)$ denote the Serre
quotient of the category of coherent $\ms C$-modules by the
subcategory consisting of modules supported in codimension $2$.  (The
peculiar notation is meant to make this align with section 1.6 of
\cite{h-l}.  The reader looking for more details about this type of category
is referred there.)

\begin{prop}\label{P:alg-stab}
  A torsion free coherent $\ms X$-twisted sheaf $\ms F$ is
  $\mu$-(semi)stable if and only if $M(\ms F)$ is a $\mu$-(semi)stable
  $\ms B$-module on $\widebar X$.
\end{prop}
\begin{proof}
  First, it is easy to see that $\ms F$ is (semi)stable if and only if
  $\sigma_{\ast}\shom(\ms V,\ms F)$ is a (semi)stable
  $\sigma_{\ast}\send(\ms V)$-module.  Thus, the content of the lemma
  is in the statement that $\pi_{\ast}$ reflects (semi)stability.

Let $\ms A$ be a locally free sheaf of $\ms O_X$-algebras with
(reflexive) pushforward $\ms B$ on $\widebar X$.  Note that $X$ and
$\widebar X$ are isomorphic in codimension $1$, which implies that the functor
between Serre quotient categories
$$\Coh_{2,1}(\ms A)\to\Coh_{2,1}(\ms B)$$
is an equivalence.  Since $\sigma_\ast$ has vanishing higher direct
images on the category of quasi-coherent sheaves, we
know that for any $\ms A$-module $\ms E$ and any integer $m$,
$\chi(X,\ms E(m))=\chi(\widebar X,\ms E(m))$.  By Lemma
\ref{L:slope-hilbert-comp}, we see that the slope of $\ms E$ equals
the slope of $\sigma_\ast\ms E$.  We conclude that the notions of semistability
coincide under pushforward, as they are defined by the same function
on the same $K$-group.
\end{proof}

\begin{prop}\label{sec:moduli-mu-semistable}
  The substack $\Tw^s_{\ms X}(r,L,\Delta)\subset\Tw^{ss}_{\ms
    X}(r,L,\Delta)$ parametrizing $\mu$-stable twisted sheaves is
  open.  Moreover, $\Tw^{ss}_{\ms X}(r,L,\Delta)$ is an open substack
  of the stack of $\ms X$-twisted sheaves.
\end{prop}
\begin{proof}
  Lemma 3.7 of \cite{simpson} shows that the substack $\ms M^{ss}_{\ms
    B}$ (respectively, $\ms M^s_{\ms B}$) of $\mu$-semistable
  (respectively $\mu$-stable) coherent $\ms B$-modules is open in $\ms
  M_{\ms B}$.
  Applying Proposition \ref{P:alg-stab}, we see that the substack of
  $\Tw_{\ms X}$ parametrizing semistable $\ms X$-twisted sheaves is
  $M^{-1}(\ms M_{\ms B}^{ss})$, and similarly for the stable loci.
\end{proof}
In sections \ref{sec:dimens-segre-loci} and
\ref{sec:asympt-prop-moduli}, we will take up (among other things) the
question of the density of $\Tw^{s}_{\ms X}(r,L,\Delta)$ in
$\Tw^{ss}_{\ms X}(r,L,\Delta)$.

\subsubsection{Pushing forward, II: comparison of Harder-Narasimhan filtrations}
\label{sec:comp-hard-naras}

Let $\ms E$ be a torsion free sheaf on $X$.  In the spirit of the last
section, we
compare properties of the Harder-Narasimhan filtration of $\ms E$ with
those of the Harder-Narasimhan filtration of $\sigma_{\ast}\ms E$.

Recall that the Hilbert polynomial gives rise to a definition of
stability for sheaves on the normal variety $\widebar X$ (the coarse
moduli space of $X$), and a similar notion of slope stability (by
paying attention to only the first two coefficients of the reduced
Hilbert polynomial).

\begin{lem}\label{sec:comp-hard-naras-1}
  We have $\mu(\ms E)=\mu(\sigma_{\ast}\ms E)$.
\end{lem}
\begin{proof}
  This follows from Lemma \ref{L:slope-hilbert-comp}.
\end{proof}

\begin{lem}\label{sec:comp-hard-naras-2}
  A sheaf $\ms E$ is slope stable if and only if $\sigma_\ast\ms E$ is
  slope stable.
\end{lem}
\begin{proof}
  This follows immediately from Lemma \ref{sec:comp-hard-naras-1} and
  the fact that slope-stability is determined in codimension $1$,
  where $X$ and $\widebar X$ are isomorphic.
\end{proof}

\begin{prop}\label{sec:comp-hard-naras-3}
  If $0=\ms E_0\subset\ms E_1\subset\cdots\subset\ms E_k=\ms E$ is the
  Harder-Narasimhan filtration of $\ms E$ (with respect to slope
  stability), then $(\sigma_\ast\ms E)_i:=\sigma_{\ast}(\ms E_i)$
  defines the Harder-Narasimhan filtration on $\sigma_\ast\ms E$ (with
  respect to slope stability).  Moreover, the slopes of the
  subquotients are the same.
\end{prop}
\begin{proof}
  This follows from the fact that $\sigma_\ast$ sends torsion free
  sheaves to torsion free sheaves 
  condition along with Lemma \ref{sec:comp-hard-naras-2} and Lemma
  \ref{sec:comp-hard-naras-1}.
\end{proof}

\subsection{Dimension estimates at points of $\Tw^{ss}_{\ms X}(r,L,\Delta)$}
\label{sec:dimension-estimates}
 Let $\ms E$ be a semistable $\ms X$-twisted sheaf of rank $r$,
 determinant $L$, and discriminant $\Delta$.  Letting $\ext^i(\ms
 E,\ms E)_0$ denote the kernel of the trace map $\ext^i(\ms E,\ms
 E)\to\H^i(X,\ms O)$, it is standard that the miniversal deformation
 of $\ms E$ has dimension $d$ satisfying
$$\dim\ext^1(\ms E,\ms E)_0-\dim\ext^2(\ms E,\ms E)_0\leq
d\leq\dim\ext^1(\ms E,\ms E)_0$$ (see, e.g., Proposition 2.A.11 of
\cite{h-l} or the classic \cite{schl}).  Thus, a lower
bound for $d$ is provided by $\chi(X,\ms O)-\chi(\ms E,\ms
E)+\dim\hom(\ms E,\ms E)_0\geq\chi(X,\ms O)-\chi(\ms E,\ms E)$ (with
equality holding when $\ms E$ is simple).  Now, we know that the
dimension of $\Tw^{ss}_{\ms X}(r,L,\Delta)$ at $\ms E$ will be at
least $d-\dim\aut(\ms E)+1\geq d-r^2+1$ (by Lemma
\ref{L:simplest-bound}), so that we conclude that
\begin{equation}\label{Eq:lower-dim-bd}
\dim_{[\ms E]}\Tw^{ss}_{\ms X}(r,L,\Delta)\geq \chi(X,\ms
O_X)-\chi(\ms E,\ms E)-r^2+1.
\end{equation}

Using the Riemann-Roch formula, we find the following.

\begin{lem}\label{sec:dimens-estim-at}
  Let $\ms E$ be a totally regular semistable $\ms X$-twisted sheaf of
  rank $r$, determinant $L$, and discriminant $\Delta$.  There is a
  lower bound
$$\dim_{[\ms E]}\Tw^{ss}_{\ms
  X}(r,L,\Delta)\geq\Delta-(r^2-1)\chi(X,\ms O_X)+r^2n\frac{N^2-1}{12N}-r^2+1.$$
\end{lem}
\begin{proof}
  This results from combining Lemma \ref{sec:discr-bogom-ineq-2} and
  equation \eqref{Eq:lower-dim-bd} above.
\end{proof}

A similar argument, using the fact that stable sheaves are simple,
shows that for any stable sheaf $\ms E$ there is an inequality

$$\Delta-(r^2-1)\chi(X,\ms O_X)+r^2n\frac{N^2-1}{12N}\leq\dim_{[\ms
  E]}\Tw^s_{\ms X}(r,L,\Delta)\leq \Delta-(r^2-1)\chi(X,\ms
O_X)+r^2n\frac{N^2-1}{12N}+\dim\ext^2(\ms E,\ms E)_0.$$

\begin{defn}\label{sec:dimens-estim-at-1}
  Given $r$, $L$, and $\Delta$, the \emph{expected
    dimension\/} of $\Tw^s_{\ms X}(r,L,\Delta)$ is
$$\exp\dim\Tw^s_{\ms X}(r,L,\Delta):=\Delta-(r^2-1)\chi(X,\ms O_X)+r^2n\frac{N^2-1}{12N}.$$
\end{defn}
We see that for a stable $\ms X$-twisted sheaf $\ms E$ we have
\begin{equation}\label{eq:2}
\exp\dim\Tw^s_{\ms X}(r,L,\Delta)\leq\dim_{[\ms E]}\Tw^s_{\ms
  X}(r,L,\Delta)\leq\exp\dim\Tw^s_{\ms X}(r,L,\Delta)+\dim\ext^2(\ms
E,\ms E)_0.
\end{equation}

\subsection{Dimensions of Segre loci}
\label{sec:dimens-segre-loci}

In this section we follow ideas established by O'Grady and Langer and
produce bounds on the set of sheaves whose Segre invariant has a fixed
upper bound.

\subsubsection{An estimate}
\label{sec:some-estimates}

We use the function $f:\Z/N\Z\to\Q$ from section
\ref{sec:correction-terms}. The sequence $f(0)$, $f(1)$,\ldots, $f(N-1)$
gives rise to a (symmetric) circulant matrix
$(c_{ij})=(\frac{1}{N}f(j-i))$.  It is well-known 
that $(c_{ij})$ is diagonalizable (over $\R$) with eigenvalues
$\eps_m=\sum_{i=0}^{N-1}\frac{1}{N}f(i)\zeta^{mi}$.

\begin{lem}
  The eigenvalues $\eps_m$ described above are all non-negative.
\end{lem}
\begin{proof}
First, a preliminary reduction.  If $m=0$, this is a calculation which
is left to the reader.  For $m>0$, the sum is taken over powers of a
root of unity $z$ such that $z^d=1$ for some $d$ dividing $N$.  The
sum simplifies to 
$$\eps_m=\sum_{i=0}^{N-1}\frac{i(i-N)}{2N}z^i=\frac{1}{2N}\left(\sum
i^2z^i-N\sum iz^i\right).$$
It thus suffices to evaluate the exponential sums $\sum i^rz^i$ for
$r=1$ and $r=2$.

We recall
  the standard power sums
$$\sum_{j=1}^{N-1}jx^j=\frac{(N-1)x^{N+1}-Nx^N+x}{(x-1)^2}$$
and
$$\sum_{j=1}^{N-1}j^2x^j=\frac{(N-1)^2x^{N+2}+(-2N^2+2N+1)x^{N+1}+N^2x^N-x^2-x}{(x-1)^3}$$
which can be derived by differentiating the standard geometric sum
$\sum_{j=1}^{N-1} x^j=(x^N-1)/(x-1)$.  Using the fact that $z^d=1$, we
find that 
$$\sum_i iz^i=\frac{N}{z-1}$$
and
$$\sum_i i^2z^i=\frac{N(N-2)z^2+2N(1-N)z}{(z-1)^3}.$$
A simple calculation yields
$$\eps_m=\frac{-2z}{(z-1)^2},$$
which we wish to show is positive.  (The dependence on $m$ is hidden
in the choice of $z$.)  It is elementary to see that the argument of
$z/(z-1)^2$ is $\pi$, which establishes the result.  (E.g., one can
easily show that it suffices to calculate the argument for any
primitive such root of unity, as the corresponding sums for different
choices of $z$ differ by positive real scalars; for $z=e^{2\pi i/d}$,
the calculation is an exercise in plane geometry, at least when $d>4$.)
\end{proof}

In other words, the matrix $(c_{ij})$ is positive semidefinite. 

Now suppose $F$ is a perfect complex on $X$.  Using the
To\"en-Riemann-Roch formula, we can write  
$$\chi(F)=\deg(\chern(F)\cdot\td_X)+\sum_{i=1}^n\delta_i(F),$$
where $\delta_i(F)$ is a correction term (a priori lying in
$\Q(\m_N)$, but actually lying in $\Q$) coming from contributions at
the residual gerbe $\xi_i$.  Write
$[\LLL\iota_i^{\ast}F]=\sum_{j=0}^{N-1}e^{(i)}_j\chi^j$ for the class in
$K$-theory of the derived fiber of $F$ over $\xi_i$.

In this section, we will be interested in computing the correction
terms to the Riemann-Roch formula for particular perfect complexes $F$ 
arising as follows: let $E$ be a totally regular twisted sheaf of rank
$r$ and
\begin{equation}\label{Eq:seq}
0\to E_1\to E\to E_2\to 0
\end{equation}
a non-trivial exact sequence of torsion free sheaves.  Let $F=\rshom(E_1,E_2)$.

\begin{prop}\label{sec:some-estimates-1}
  Using the above notation, for any exact sequence \eqref{Eq:seq} and any 
  $i=1,\ldots,n$ we have that $\delta_i(F)\leq 0$.
\end{prop}
\begin{proof}
  We work at a fixed residual gerbe $\xi_i$ and omit $i$ from the
  notation.  Using the notation of section
  \ref{sec:sheav-deform-theory}, for a sheaf $G$ we will write $\widebar{G}$ for
  $\widebar G_i$.  Let $\widebar{E_1}=\sum d_i\chi^i=:P$, and let $\bf
  d$ be the column vector $(d_i)$.  
Let $s=\sum d_i<N$ be the rank of $E_1$.  We
  claim that
$$\delta(F)=-\mathbf d^{\operatorname T}(c_{ij})\mathbf d+\frac{s}{N}\sum_{j=0}^{N-1}f(j).$$
Then, since $(c_{ij})$ is positive semidefinite, it follows that
$\delta(F)$ is bounded by $(s/N)\sum f(j)=0$.

To prove the claim, note that there is an equality $\widebar{E_2}=\rho-P$ in
$K(\m_N)$, so that $\widebar F=P^{\vee}\tensor
(\rho-P)=r\rho-P^{\vee}\tensor P$.  We can write $P^{\vee}\tensor
P=\sum_{i,j}d_id_{i+j}\chi^j$.  Thus, since $\chi(F)$ depends only on
the class of $F$ in $K$-theory, the correction term becomes
$$\delta(F)=\frac{s}{N}\sum
f(j)-\sum_{j=0}^{N-1}\left(\frac{1}{N}f(j)\sum_{i=0}^{N-1}d_id_{i+j}\right).$$
The second sum in this expression is a quadratic form in the $d_i$ in which the
coefficient of $d_ad_b$ is given by $-\frac{1}{N}f(b-a)-\frac{1}{N}f(a-b)=-\frac{2}{N}f(b-a)$ when $a\neq b$ and
by $-\frac{N^2-1}{12N}=-f(0)$ when $a=b$.  The matrix corresponding to this form is precisely $-(c_{ij})$,
as desired.
\end{proof}

\subsubsection{Flag spaces of semistable twisted sheaves}
\label{sec:flag-spac-semist}

We recall a few basic facts on dimensions of substacks of the stack of
semistable twisted sheaves which admit a specified flag structure.  We have
followed paragraph 3.7 of \cite{langer-castelnuovo} and appendix 2.A
of \cite{h-l} in our presentation, with a few
modifications appropriate to the situation at hand.

Given rational polynomials $P_1,\ldots,P_k$, let
$\Filt(P_1,\ldots,P_k)$ denote the stack of filtered $\ms X$-twisted
sheaves $F_0\ms E\subset F_1\ms E\subset\cdots\subset F_k\ms E=\ms E$
such that the $f$-Hilbert polynomial of $\ms E_i:=F_{i}\ms
E/F_{i-1}\ms E$ is $P_i$ for $i=1,\ldots,k$.  The map $F_i\ms
E\mapsto\ms E/F_0\ms E$ defines a morphism of Artin stacks 
$$q:\Filt(P_1,\ldots,P_k)\to\Tw_{\ms X}.$$

\begin{lem}
  The reduced structure on the image of $q$ is a closed substack of $\Tw_{\ms X}$.
\end{lem}
\begin{proof}
  It is enough to check this after pulling back to a smooth cover of
  $\Tw_{\ms X}$, so this reduces to the following: if $\ms E$ is an
  $S$-flat coherent $\ms X\times S$-twisted sheaf with torsion free
  fibers, then the locus parametrizing sheaves admitting a flag of the
  prescribed type is locally closed in $S$.  We proceed by induction
  on $k$, the case when $k=0$ following from the fact that the
  geometric Hilbert polynomial is constant in a flat family.  We know
  that $Q:=\Quot(\ms E,P_k)$ is proper over $S$ by Corollary
  \ref{sec:unif-m_n-gerb}.  Over $Q$, there is a two-step filtration $F_{k-1}\ms
  E_Q\subset\ms E_Q$ with quotient having geometric Hilbert polynomial
  $P_k$.  By induction on $k$, the subset of $Q$ containing the rest
  of the flag is closed; since $Q$ is proper over $S$, the image is
  closed, as desired.
\end{proof}

We will (temporarily) denote the image of $q$ by $\Tw_{\ms X}(P_1,\ldots,P_k)$
(following Langer).

\begin{prop}\label{sec:flag-spac-semist-1}
  The dimension of $\Tw_{\ms X}(P_1,\ldots,P_k)$ at a point $\ms E$
  admitting a filtration $0=F_0\ms E\subset F_1\ms
  E\subset\cdots\subset F_k\ms E=\ms E$ is bounded above by
$$\sum_{i\geq j}\dim\ext^1(\ms E_i,\ms E_j)+\dim\aut(\ms E)-1,$$
where $\ms E_i:=F_i\ms E/F_{i-1}\ms E$.  If $\ms E$ is semistable of
rank $r$ then $\dim\aut(\ms E)\leq r^2$.
\end{prop}
\begin{proof}
  Fix a locally free $\ms X$-twisted sheaf $\ms V$.  We can also
  realize the stack $\Tw_{\ms X}$ as the stack
  $\operatorname{Mod}_{\ms A}$ of coherent right 
  $\ms A$-modules (by sending $\ms E$ to $\pi_{\ast}\shom(\ms V,\ms
  E)$, where $\ms A=\pi_{\ast}\send(\ms V)$ is an Azumaya algebra on
  $X$.  Given an $\ms A$-module $\ms F$, we also have the usual
  Hilbert polynomial $P_{\ms F}(m)=\chi(X,\ms F(m))$, where we twist
  by a fixed polarization of the coarse moduli space of $X$; as usual,
  for large values of $m$ we have $P_{\ms F}(m)=\H^0(X,\ms F(m))$ (as
  Serre's theorem applies).

  The stack $\Tw_{\ms X}$ is a union of open substacks $\Tw_{\ms
    X}^{(m)}$ consisting of twisted sheaves $\ms E$ such that
  $\shom(\ms V,\ms E)$ is $m$-regular in the sense of Mumford.  We can
  realize any such $\ms E$ as a quotient of $$\hom(\ms V,\ms
  E(m))\tensor_k\ms V(-m).$$ Thus, if we fix the Hilbert polynomial
  $P_{\shom(\ms V,\ms E)}$ we get a quasi-compact stack.  It follows
  that the stack $\Filt(P_1,\ldots,P_k)^{(m)}$ of filtrations on
  $m$-regular twisted sheaves breaks up as a finite disjoint union of
  corresponding stacks of filtered $\ms A$-modules (over a finite
  collection of sequences of Hilbert polynomials for the subquotient
  $\ms A$-modules), each of which is open and closed in
  $\Filt(P_1,\ldots,P_k)$.

  Thus, to estimate the dimension of $\Tw_{\ms X}(P_1,\ldots,P_k)$, we
  can also assume that the filtrations in question have fixed $\ms
  A$-Hilbert polynomials $Q_1,\ldots,Q_k$, and that the ambient
  sheaves have Castelnuovo-Mumford regularity at most $m$.  Let $\ms
  H=\ms A(-m)^{Q_{\shom(\ms V,\ms E)}(m)}$.  Standard methods (just as
  in Proposition 3.8 of \cite{langer-castelnuovo}) show that the flag
  space $Y$ parametrizing $\ms A$-module filtrations $G_0\ms
  H\subset\cdots\subset G_k\ms H=\ms H$ such that the Hilbert
  polynomials of the subquotients are $P_i$ and the $\ms A$-Hilbert
  polynomials are $Q_i$ has dimension at most $$\sum_{i\geq
    j}\dim\ext^1(\ms E_i,\ms E_j)+\dim\End(\ms H)-1.$$ On the other
  hand, there is a map $\alpha:Y\to\Tw_{\ms X}$ sending a filtration
  to $\ms H/G_0\ms H$. This morphism obviously factors through the map
  $\beta:Y\to Q:=\Quot(\ms V(-m)^{Q_{\shom(\,s V,\ms E)(m)}},P)$; moreover,
    the image of $\beta$ is the preimage of the image of $\alpha$.
    Thus, to get an upper bound on the dimension of $\Tw_{\ms
      X}(P_1,\ldots,P_k)$ at $\ms E$, it suffices to subtract from
    $\dim Y$ a lower bound for the dimension of the fibers of
    $Q\to\Tw_{\ms X}$ near $\ms E$.  Elementary
    considerations show that the fiber dimension at $\ms E$ is given
    by $\dim\End_{\ms A}(\ms H)-\dim\Aut(\ms E)$, which gives the desired
    result.

  If $\ms E$ is semistable, then $\dim\hom(\ms E,\ms E)\leq r^2$ by
  \ref{L:simplest-bound}, from which it follows that $\dim\aut(\ms E)\leq r^2$.
\end{proof}

\subsubsection{The dimension estimate}
\label{sec:dimension-estimate}

This section is again very close to the approach taken by Langer in
section 7 of \cite{langer-castelnuovo}.  We recall the principal
results and give proofs only when they are different from those of
Langer.  Following Langer, define $f(r)=-1+\sum_{i=1}^r\frac{1}{i}$.

We retain the notation from section \ref{sec:discr-bogom-ineq}.  Given
a torsion free coherent $\ms X$-twisted sheaf $F$, it follows from Lemma \ref{L:i-cover-you} that the Harder-Narasimhan filtration on $F$ pulls back to the
Harder-Narasimhan filtration on $f^\ast F$.  The calculations of
section \ref{sec:discr-bogom-ineq} show that
$\mu_{max}(f^{\ast}F)=d\mu_{\max}(F)$.

\begin{lem}\label{sec:dimension-estimate-2}
  Let $F_1$ and $F_2$ be torsion free slope-semistable coherent $\ms
  X$-twisted sheaves such that the rank and slope of $F_i$ are $r_i$
  and $\mu_i$, respectively.  There is an inequality
$$\mu_{max}(\shom(F_1,F_2))\leq\mu_2-\mu_1+(r_1+r_2-2)L_X.$$
\end{lem}
\begin{proof}
  It is clear that we can replace $F_1$ and $F_2$ by their reflexive
  hulls and thus assume that they are locally free.  Thus, we may
  assume that the formation of $\shom(F_1,F_2)$ commutes with pullback
  to $Y$. Taking into account the fact that $L_Y=dL_X$ and similarly
  for the slopes in question, the fact to be proved follows follows immediately
  from the corresponding fact for sheaves on smooth projective
  surfaces.  This is contained in the proof of Corollary 4.2 of
  \cite{langer-castelnuovo} (and is a short deduction from equation (2.1.1) of
  [\emph{ibid\/}.]).
\end{proof}

\begin{lem}\label{sec:dimension-estimate-1}
  Let $E_1$ and $E_2$ be torsion free slope-semistable $\ms X$-twisted
  sheaves.  Let $r_i$ and $\mu_i$ denote the rank and slope of $E_i$.
  If $\mu_1>\mu_2$ then $\hom(E_1,E_2)=0$.  Otherwise, we have
$$\dim\hom(E_1,E_2)\leq\frac{r_1r_2}{2H^2}\left(\mu_2-\mu_1+\frac{3}{2}H^2+(r_1+r_2-2)L_{X}\right)^2+\frac{r_1r_2}{8}(2f(r_1r_2)-9)H^2.$$
\end{lem}
\begin{proof}
  The proof is identical to the proof of Corollary 4.2 in
  \cite{langer-castelnuovo}, using Lemma \ref{sec:dimension-estimate-2} to bound
  $$\mu_{max}(\shom(E_1,E_2))=\mu_{max}(\sigma_{\ast}\shom(E_1,E_2))$$
  and Theorem 4.1 of [\emph{ibid\/}.] to bound $\H^0(\widebar
  X,\sigma_\ast\pi_\ast\shom(E_1,E_2))$.
\end{proof}

\begin{cor}\label{sec:dimension-estimate-5}
  There exists a constant $\beta_{\infty}$ depending only on $X$ and
  $r$ such that for any semistable $\ms X$-twisted sheaf $\ms F$ of
  rank $r$ invertible in $k$ we have 
$$\dim\ext^2(\ms F,\ms F)_0\leq\beta_\infty.$$
\end{cor}
\begin{proof}
  This follows immediately from Lemma \ref{sec:dimension-estimate-1}
  with $E_1=\ms F$ and $E_2=\ms F\tensor\omega_X$.
\end{proof}

\begin{cor}\label{sec:dimension-estimate-6}
  There exists a constant $\beta(r)$ depending only on $X$ and $r$ such that
  for any semistable $\ms X$-twisted sheaf $\ms E$ of rank $r$ and determinant
  $L$, we have
$$\dim_{\ms E}\Tw^{ss}_{\ms X}(r,L,\Delta)\leq\Delta+\beta(r).$$
\end{cor}
\begin{proof}
  The proof is identical to the proof of Lemma 8.5 of \cite{langer-castelnuovo}.
\end{proof}

\begin{cor}\label{sec:dimension-estimate-7}
  We have $\dim\Tw^s_{\ms X}(r,L,\Delta)\leq\exp\dim\Tw^s_{\ms X}(r,L,\Delta)+\beta_\infty$.
\end{cor}
\begin{proof}
  This follows from Corollary \ref{sec:dimension-estimate-5} and
  equation (\ref{eq:2}).
\end{proof}

Using \ref{sec:dimension-estimate-1}, we can now give an estimate of
the dimension of the locus $\Tw_{\ms X}^{ss}(r,L,\Delta)(s)$ parametrizing
semistable $\ms X$-twisted sheaves with Segre invariant at most $s$.

\begin{lem}\label{sec:dimension-estimate-4}
  Let $E$ be a slope semistable $\ms X$-twisted sheaf of rank $r$.
  Let the Segre invariant $s(E)$ be realized by an exact sequence
$$0\to E_1\to E\to E_2\to 0,$$
and let $r_i:=\rk E_i$.  There exists a constant $\alpha(r)$ depending
only on $X,\ms X,Y,H$ and $r$ such that
$$\sum_{i\geq j}\dim\ext^1(E_i,E_j)\leq\frac{\max(r_1,r_2)+r}{2r}\left(\Delta(E)+\frac{r_1r_2}{H^2}s^2\right)+\frac{r_1r_2}{2H^2}(3H^2-KH+2(r-2)L_X)s+\alpha(r).$$
\end{lem}
\begin{proof}
  The reader should compare this with Proposition 7.1 of
  \cite{langer-castelnuovo}.  The proof is identical, with the following
  modifications:  First, in place of Langer's Bogomolov inequality
  (Theorem 2.2 of [\emph{ibid\/}.]) we
  use Proposition \ref{sec:discr-bogom-ineq-1} (and the orbifold Hodge
  Index Theorem \ref{T:pre-hodge} above).  In place of Langer's
  Corollary 4.2 of [\emph{ibid\/}.], we use Lemma
  \ref{sec:dimension-estimate-1}.  The final modification requires
  that we explain the main outline of the proof.

We start with the equality
$$\sum_{i\geq j}\dim\ext^1(E_i,E_j)=\sum_{i\geq
  j}\left(\dim\hom(E_i,E_j)+\dim\hom(E_j,E_i\tensor\omega_X)\right)+\chi(E_1,E_2)-\chi(E,E).$$
Lemma \ref{sec:dimension-estimate-1} yields bounds for the terms in
the sum on the right-hand side, so it remains to bound the difference
$\chi(E_1,E_2)-\chi(E,E)$, and this is the point where one must modify
the proof.  The following lemma plays a central role.

\begin{lem}
  Given a perfect complex $F$ of coherent $\ms O_X$-modules of
  positive rank $r$, 
we have
$$\chi(F)=-\frac{\Delta(F)}{2r}+\frac{r}{2}\xi(\xi-K_X)+\chi(\ms O_X)+c(F),$$
where $\xi=c_1(F)/r$ and $c(F)$ is a correction term depending only upon the classes
$[\LLL\iota_i^{\ast}F]$ in $K(\B\m_N)$.  If
$F=\R\sigma_{\ast}\shom(E_1,E_2)$ with $0\to E_1\to E\to E_2\to 0$ an
exact sequence of twisted sheaves as in \eqref{Eq:seq} above, then
$c(F)$ is bounded in terms of $X$,$N$,$n$, and $r$.
\end{lem} 
\begin{proof}
  The To\"en-Riemann-Roch formula shows that
  $\chi(F)=\deg(\chern(F)\cdot\td_X)+c'(F)$, where $c'(F)$ is a
  correction term depending only upon the $[\LLL\iota_i^{\ast}F]$.
  Expanding out the intersection product in Chern classes and using
  Corollary \ref{C:euler-structure-sheaf} yields the first part of the
  lemma.  The second follows from Proposition
  \ref{sec:some-estimates-1}.
\end{proof}

Let $\xi=\widetilde c_1(E_1)/r_1-\widetilde c_1(E_2)/r_2$.  The
preceding lemma shows that 
$$\chi(E_1,E_2)-\chi(E,E)=\Delta(E)-\left(r_2\frac{\Delta_1}{2r_1}+r_1\frac{\Delta_2}{2r_2}\right)+\frac{1}{2}r_1r_2(\xi-K_X)\xi+(r_1r_2-r^2)\chi(\ms
O_X)+c,$$ where $c$ is a correction term which is bounded above in
terms of $X$,$N$,$n$, and $r$.  The proof of Proposition 7.1 of
\cite{langer-castelnuovo} yields an upper bound on
$\chi(E_1,E_2)-\chi(E,E)-c$ (by following the arguments verbatim);
since $c$ is bounded, this yields a bound on $\chi(E_1,E_2)-\chi(E,E)$
(and changes Langer's constant $\alpha_2(r)$ by a constant).  The rest
is precisely as in the proof of [\emph{loc.\ cit\/}.].
\end{proof}

\begin{prop}\label{sec:dimension-estimate-3}
  There exists a constant $B$ such that for any $s\geq 0$
$$\dim\Tw_{\ms X}^{ss}(r,L,\Delta)(s)\leq\left(1-\frac{1}{2r}\right)\Delta+\frac{2r^2}{9H^2}s^2+\frac{r^2}{8H^2}\left[3H^2-KH+2(r-2)L_X\right]_{+}s+B.$$
\end{prop}
\begin{proof}
  The proof is identical to the proof of Theorem 7.2 in
  \cite{langer-castelnuovo}, using Proposition
  \ref{sec:flag-spac-semist-1} in place of Langer's Proposition 3.8
  and Lemma \ref{sec:dimension-estimate-4} in place of Langer's
  Proposition 7.1
\end{proof}

\section{Asymptotic properties of moduli}
\label{sec:asympt-prop-moduli}
\setcounter{subsubsection}{0}

In this section we study what happens to the stack $\Tw^s_{\ms
  X}(r,L,\Delta)$ as $\Delta$ goes to infinity.

\subsection{Non-emptiness}
\label{sec:non-emptiness}

\begin{defn}\label{sec:non-emptiness-3}
  Given a torsion free $\ms X$-twisted sheaf $\ms V$, the
  \emph{discrepancy\/} of $\ms V$ is 
$$\delta(\ms V):=\max_{\ms V'\subset\ms V}\{\mu(\ms V')-\mu(\ms
V)\},$$
where $\ms V'\subset\ms V$ ranges over proper subsheaves satisfying
$0<\rk\ms V'<\rk\ms V$.
\end{defn}
It is clear that $\delta(\ms V)$ has a denominator bounded by
$N\rk(\ms V)!$.  If $\ms V$ is unstable then $\delta(\ms
V)=\mu_{\text{\rm max}}(\ms V)-\mu(\ms V)$, while $\ms V$ is stable
if and only if $\delta(\ms V)<0$.  It is elementary that the discrepancy is
invariant under twisting by an invertible sheaf.

\begin{prop}\label{sec:non-emptiness-1}
  Suppose $\ms V$ is a torsion free totally regular $\ms X$-twisted
  sheaf of rank $r$ with determinant $L$.  Given any $c\in\Q$, there exists
  $c>c_0$ and a subsheaf $\ms W\subset\ms V$ with rank
  $r$, determinant $L$ and $\deg c_2(\ms V)>c$ such that $\delta(\ms
  W)<\delta(\ms V)$.  In particular, there is such a subsheaf $\ms W$
  which is $\mu$-stable.
\end{prop}
\begin{proof}
  This is similar to Theorem 5.2.5 of \cite{h-l}.  The reader
  will note that the proof of [\emph{loc.\ cit.\/}] uses the fact that
  the trivial bundle is semistable, a fact which is not available to
  us.  The proof here shows that even if we start with an unstable
  ambient sheaf $\ms V$, we can at least ensure that the subsheaf will
  have smaller discrepancy.  Since the discrepancies have bounded
  denominator, iterating this process will ultimately produce a stable sheaf.

  A preliminary reduction: replacing $\ms V$ by $\ms V(-tH)$ for large
  $t$, we may assume that $$\mu(\ms V)<\min\{0,(r-1)\delta(\ms V)\}.$$  Since
  the discrepancy of a sheaf is invariant under twisting, we can twist
  the resulting subsheaf of $\ms V(-tH)$ back up by $tH$ to achieve
  the desired result.

  Let $C$ be a smooth curve belonging to the linear system $|rmH|$ for
  large $m$ which lies in $X^\sp$ and avoids the singular points of
  $\ms V$; write $\ms C:=\ms X\times_X C$.
  A result of Grothendieck (Lemma 1.7.9 of \cite{h-l}) shows that the
  space $Q$ of locally free quotients $F$ of $\ms V|_{\ms C}$ with
  $\rk(F)>r$ and $\mu(F)\leq\frac{r-1}{r}C\cdot H$ is quasi-compact.
  Let $\ms M$ be an invertible $\ms C$-twisted sheaf (which exists by
  Tsen's theorem). Let $q:\pr_1^{\ast}\ms V_{\ms C}\surj\ms F$ on $\ms
  C\times Q$ be the universal quotient of $\ms V|_{\ms C}$ satisfying
  those two conditions.  For large $s$ we know that $(\pr_2)_\ast\left(\ms
  F^{\vee}\tensor\pr_1^{\ast}\ms M(sH)\right)$ is a locally free sheaf on
  $Q$; it is thus the sheaf of sections of some geometric vector
  bundle $\V\to Q$ whose fiber over a quotient $\ms V|_{\ms C}\to F$
  is $\hom(F,\ms M(sH))$.  The quotient map $q$ induces a map
  $\widetilde q:\V\to\uhom(\ms V|_{\ms C},\ms M(sH))$ which is a
  linear embedding when restricted to each fiber over $Q$.  Since
  $\rk(\ms F)<r$, we see that for sufficiently large $s$ the map
  $\widetilde q$ cannot be dominant.  Thus, there exists some
  surjective map
  $\phi:\ms V\to\ms M(sH)$ which does not factor through any
  quotient parametrized by $Q$.

  We claim that the kernel $\ms W$ of one such map $\phi:\ms V\to\ms
  M(sH)$ satisfies $\delta(\ms W)<\delta(\ms V)$.  Let us grant this
  for a moment.  Standard computations show that $\det(\ms W)=\det(\ms
  V)(-C)$, so that the sheaf $\ms W(-mH)$ will then have determinant
  $L$ and $\delta(\ms W)=\delta(\ms W(-mH)$.  Taking the kernel of a
  general quotient $\ms W\to\ms Q$ with $\ms Q$ of finite length will
  then make the second Chern class arbitrarily large.

  Thus, it remains to establish the assertion about the disrepancy.
  Given a saturated subsheaf $\ms W'\subset\ms W$, let $\ms
  V'\subset\ms V$ be its saturation in $\ms V$.  The sheaf $\ms R:=\ms
  V'/\ms W'$ is a subsheaf of $\ms M(sH)$.  If $\ms R\neq 0$ then we
  have that $\det(\ms W')\cong\det(\ms V')(-C)$, so that
$$\mu(\ms W')=\mu(\ms V')-\frac{C\cdot H}{\rk(\ms W')}<\delta(\ms
V)+\mu(\ms V)-\frac{C\cdot H}{r}=\delta(\ms V)+\mu(\ms W),$$ so that
$\mu(\ms W')-\mu(\ms W)<\delta(\ms V)$.  If $\ms R=0$, then the sheaf $\ms
U=\ms V/\ms W'$ is a torsion free $\ms C$-twisted sheaf, and the map
$\phi:\ms V\to\ms M(sH)$ factors through the quotient $\ms U$.  By
construction, we know that $\mu(\ms U)\geq\frac{r-1}{r}C\cdot H$.
From the exact sequence
$$0\to\ms W'\to\ms V\to\ms U\to 0$$ we deduce that
\begin{align*}
\mu(\ms W')&=\mu(\ms V)-\frac{r-\rk\ms W'}{\rk\ms W'}\mu(\ms U)\\
&\leq\mu(\ms V)-\frac{r-\rk\ms W'}{\rk\ms W'}\frac{r-1}{r}C\cdot H\\
&<\frac{1}{r-1}\mu(\ms V)+\mu(\ms W),
\end{align*}
from which it follows that $\mu(\ms W')-\mu(\ms W)<\delta(\ms V)$, as desired.
\end{proof}

\begin{cor}\label{sec:non-emptiness-2}
  Given $r$, $L$, and $c_0$, if there exists a totally regular torsion
  free $\ms X$-twisted sheaf $\ms V$ of rank $r$ and determinant $L$
  then $\Tw^s_{\ms X/k}(r,L,c)\neq\emptyset$ for some $c>c_0$.  If in
  addition $\ms V$ is locally free at the non-trivial residual gerbes
  of $X$ then there is a point of $\Tw^s_{\ms X/k}(r,L,c)$ for some
  $c>c_0$ parametrizing a sheaf which is locally free at the
  non-trivial residual gerbes of $X$.
\end{cor}
\begin{proof}
  The first statement follows directly from Proposition
  \ref{sec:non-emptiness-1}.  The second follows from the proof:
  modifications are made along general curves and points, and thus do
  not change the sheaf on the (discrete) stacky locus of $X$.
\end{proof}

\subsection{Shrinking of the boundary}
\label{sec:asymp-shrink-bound}

Given a locally closed substack $Z\subset\Tw^{ss}_{\ms
  X}(r,L,\Delta)$, we will write $Z^s$ for the open substack
$Z\cap\Tw^s_{\ms X}(r,L,\Delta)$.

\begin{prop}
  There exists a constant $D$ such that for $\Delta\geq D$, the
  open substack $\Tw^{s}_{\ms X}(r,L,\Delta)\subset\Tw^{ss}_{\ms
    X}(r,L,\Delta)$ is everywhere dense.
\end{prop}
\begin{proof}
  By definition, the complement $\Tw^{ss}_{\ms
    X}(r,L,\Delta)\setminus\Tw^{s}_{\ms X}(r,L,\Delta)$ is precisely
  $\Tw^{ss}_{\ms X}(r,L,\Delta)(0)$ (the locus parametrizing
  semistable twisted sheaves with nonpositive Segre invariant).  By
  Proposition \ref{sec:dimension-estimate-3} we have that 
$$\dim\Tw^{ss}_{\ms
  X}(r,L,\Delta)\leq\left(1-\frac{1}{2r}\right)\Delta+B.$$
On the other hand, we know from Lemma \ref{sec:dimens-estim-at} that
$\dim\Tw^{ss}_{\ms X}(r,L,\Delta)\geq\Delta-(r^2-1)\chi(X,\ms
O_X)+r^2n\frac{N^2-1}{12N}-r^2+1$.  When $\Delta$ is sufficiently large the latter is
strictly larger than the former, which implies the result.
\end{proof}

\begin{hyp}\label{sec:asympt-prop-moduli-1}
  For the rest of this section, we will assume that $\Delta\geq D$, so
that every irreducible component of $\Tw^{ss}_{\ms X}(r,L,\Delta)$
contains geometrically $\mu$-stable twisted sheaves.
\end{hyp}

Fix a positive integer $n$ and let $C\subset X^\sp$ be a general curve
in $|\ms O_X(nH)|$.  Write $\ms C=\ms X\times_X C$.

\begin{prop}\label{P:restr-issue}
    If $Z$ is a closed irreducible substack of
  $\Tw^{ss}_{\ms X/k}(r,L,\Delta)$ which meets the stable
  locus such that ${\softboundary}Z=\emptyset$ and $\dim
  Z>\dim\Tw^{ss}_{\ms C}(r,L_C)$ then there is a point of
  $Z$ whose restriction to $\ms C$ is not stable.
\end{prop}
\begin{proof}
  The proof is exactly same as the proof of Lemma 3.2.4.13 of
  \cite{twisted-moduli}.  Let us recapitulate the main points.

  Suppose that every sheaf in $Z$ has (geometrically) stable
  restriction to $C$, so that there is an induced map
  $Z\to\Tw^s_{\ms C}(r,L_C)$.  Let $\widebar Z$ denote the sheafification
  of $Z^s$ (so that $Z^s\to\widebar Z$ is a $\G_m$-gerbe) and let
  $\mTw$ denote the coarse moduli space of $\Tw^s_{\ms C}(r,L_C)$ (so that,
  again, we have a gerbe).  There is a diagram
$$\xymatrix{Z\ar[r] & \Tw^s_{\ms C}(r,L_C)\ar[dd]\\
Z^s\ar[u]\ar[d] & \\
\widebar Z\ar[r]^{f}& \mTw.}$$
By assumption, $\dim\widebar Z>\dim\mTw$, so that a general geometric
fiber $F$ of $\widebar Z\to\mTw$ has dimension at least $1$.  Let
$D^\circ$ denote the (possibly affine) normalization of a connected
component of $F$.

By Tsen's theorem, there is a generically (on $D^\circ$) separable map
$\phi:D^\circ\to Z^s$ such that the induced map $D^\circ\to\Tw^s_{\ms
  C}(r,L_C)$ is essentially constant.  Using Langton's theorem
(appendix 2.B of \cite{h-l}, which easily adapts to the situation at
hand and magically does not require any extension of the fraction
field to achieve a semistable extension), we can extend $\phi$ to a
map $\phi':D\to Z$, where $D$ is the smooth projective completion of
$D^\circ$.  (Indeed, Langton's theorem extends the family corresponding to
$\phi$ to a family of semistable twisted sheaves over the base $D$;
the fact that $Z$ is closed in $\Tw^{ss}_{\ms X}(r,L,\Delta)$ then
shows that the limits along $D$ must lie in $Z$.  This is true even
though the ambient stack $\Tw^{ss}_{\ms X}(r,L,\Delta)$ is highly
non-separated in general.)  Since $\mTw$ is separated, we see that
$\phi'$ collapses $D$ to a point of $\mTw$, giving rise to an
essentially constant family of restrictions.

Thus, we have a family $\ms F$ on $\ms X\times D$ such that the
restriction to $\ms C\times D$ has fibers which are all mutually
isomorphic stable locally free $\ms C$-twisted sheaves.  By an
analysis identical to that of Lemma 3.1.4.6ff of
\cite{twisted-moduli}, we deduce that for any positive integer $n$
there is a snc divisor $C^{(n)}$ in $|nC|$ such that the restriction
of $\ms F$ to $(C^{(n)}\times_X\ms X)\times D$ is an essentially
constant family of simple locally free sheaves.  If $\ms G$ is one of
the fibers of $\ms F$, we thus see that the induced Kodaira-Spencer
map $\tau:\H^1(X,\send(\ms G))\to\H^1(C^{(n)},\send(\ms
G|_{C^{(n)}}))$ sends the element $\eps$ corresponding to the family
parametrized by $\phi'$ to $0$.  On the other hand, since $C^{(n)}$ is
arbitrarily ample, we have that $\tau$ is injective.  It follows that
$\eps=0$; applying this argument at the points of $D$ shows that
$\phi'$ is essentially constant (as it is generically separable and
has trivial tangent map), which is a contradiction.
\end{proof}

Keeping the notation of Proposition \ref{P:restr-issue}, let $[\ms
F]\in Z$ be a sheaf with non-stable restriction to $C$, so that there is
an exact sequence
$$0\to\ms F'\to\ms F|_C\to\ms F''\to 0$$
of locally free twisted sheaves on $C$ with $\mu(\ms F')\geq\mu(\ms
F'')$.  Write $r'=\rk\ms F'$ and $r''=\rk\ms F''$.  

\begin{prop}\label{P:segre-bound}
  With the above notation, if 
$$\dim Z>\dim\Tw^{ss}_{\ms X/k}(r,L,c)+r'r''(1+C(K-C))$$
then $s(\ms F)\leq 2CH$.
\end{prop}
\begin{proof}
  The proof is identical to the proof of Proposition 8.4 of
  \cite{langer-castelnuovo}, using Lemma \ref{L:slop-stab-segre-2} in
  place of Langer's Lemma 1.4.  The reader nervous about stacky issues
  need only note that $C$ lies entirely in $X^\sp$, and Corollary
  2.2.9 of \cite{h-l} applies equally to twisted sheaves on curves.
\end{proof}

Let $\theta(r)=\frac{73}{36}r^2+r-1$.

\begin{thm}\label{T:bdry-exists}
  There exist constants $A_1,C_1,C_2\in\Q$ depending only on $X$ and
  $r$ such that if $\Delta\geq A_1$ and $Z\subset\Tw^{ss}_{\ms
    X/k}(r,L,\Delta)$ is an irreducible closed substack with 
$$\dim Z\geq \left
  (1-\frac{r-1}{2r\theta(r)}\right)\Delta+C_1\sqrt{\Delta}+C_2$$
then ${\softboundary}Z^{s}\neq\emptyset$.
\end{thm}
\begin{proof}
  The proof is identical to the proof of Theorem 8.1 of Langer, using
  Proposition \ref{P:restr-issue} and Proposition \ref{P:segre-bound}
  in place of Langer's Propositions 8.2 and 8.4, and using Corollary
  \ref{sec:dimension-estimate-6} in place of Langer's Lemma 8.5.

  There is one point which deserves a bit of elaboration.  Langer (in
  the proof of Theorem 8.1) and Huybrechts-Lehn (in the proof of their
  analogous Theorem 9.2.2) reduce to proving that
  ${\softboundary}Z\neq\emptyset$ by invoking the fact that
  $\Tw^{ss}_{\ms X}(r,L,\Delta)$ is a quotient stack.  Let us rapidly sketch
  the intrinsic form of this argument: we know that the codimension of
  the dominant component of ${\softboundary}Z$ (when it is non-empty,
  as the proof shows) has codimension at most $r-1$ by Lemma
  \ref{L:bdry-codim}.  On the other hand, we know that the dimension
  of $\Tw^{ss}_{\ms X}(r,L,\Delta)(0)$ is asympotically equal to
  $\left(1-\frac{1}{2r}\right)\Delta$, which grows more slowly than
  $\left (1-\frac{r-1}{2r\theta(r)}\right)\Delta$.  Thus, for $\Delta$
  sufficiently large we will have that the codimension of
  $\Tw^{ss}_{\ms X}(r,L,\Delta)(0)$ is larger than $\dim Z-r+1$,
  implying that the dominant component of ${\softboundary}Z$ must meet
  $\Tw^{s}_{\ms X}(r,L,\Delta)$ and therefore
  ${\softboundary}Z^{s}\neq\emptyset$.  (This is obviously related to
  the other proof by looking at any smooth cover of the stack
  $\Tw^{ss}_{\ms X}(r,L,\Delta)$.)
\end{proof}

Let
$C_3=\max\{C_2+\beta_\infty,\frac{A_1}{2r}+2\beta_\infty-(r^2-1)\chi(X,\ms
O_X)+r^2n\frac{N^2-1}{12N}\}$.  In
particular, $C_3\geq C_2$.  Given a substack $Z\subset\Tw^{ss}_{\ms
  X}(r,L,\Delta)$, let
$\beta(Z)=\min\{\dim\hom_0(E,E\tensor\omega_X):[E]\in Z\}$.  We always
have that $\beta(Z)\leq \beta_\infty$ defined in Corollary \ref{sec:dimension-estimate-5}.

Following Friedman and \cite{h-l}, we make the following definition.
\begin{defn}
  An $\ms X$-twisted sheaf $\ms F$ is \emph{good\/} if it is
  $\mu$-stable and $\ext^2(\ms F,\ms F)_0=0$.
\end{defn}

Let $W\subset\Tw^{ss}_{\ms X}(r,L,\Delta)$ be the reduced closed
substack parametrizing semistable $\ms X$-twisted sheaves $\ms F$ such
that $s(\ms F)=0$ or $\ext^2(\ms F,\ms F)_0\neq 0$.  In other words,
$W$ is the complement of the good locus.

\begin{thm}\label{T:good-big}
  For all $\Delta\geq A_1$
$$\dim
W\leq\left(1-\frac{r-1}{2\theta(r)}\right)\Delta+C_1\sqrt{\Delta}+C_3.$$
\end{thm}
\begin{proof}
  The proof again follows a standard outline, with some slight
  complications arising from the fact that we only consider the soft
  boundaries.  

  Suppose $Z$ is an irreducible component of $W$ of maximal dimension.
  Arguing by contradiction, we assume that $\dim
  Z>\left(1-\frac{r-1}{2\theta(r)}\right)\Delta+C_1\sqrt{\Delta}+C_3$.
  Applying Theorem \ref{T:bdry-exists}, we know that $\softboundary
  Z^s\neq\emptyset$.  Let $D$ be the dominant component of
  $\softboundary Z^s$, and let $\ell$ be the colength of the sheaves
  parametrized by the points of $D$.  Let $Y$ be an irreducible
  component of $D$; the formation of the soft hull yields a morphism
  $p:Y\to\Tw^{s}_{\ms X}(r,L,\Delta-2r\ell)$.  Let $Z_1$ be the closure
  of $p(Y)$.  As long as $\softboundary Z_1^s\neq\emptyset$, we can
  iterate this procedure and produce $Z_2$, etc.  Note that at each
  step, the discriminant $\Delta$ decreases.  By Proposition
  \ref{sec:discr-bogom-ineq-1}, this process must terminate (as the
  discriminant is bounded below), say at $Z_t$.  

We can also understand this process by ``running it in reverse.''
That is to say: a general point $\ms E$ of $Z_1$ is softly reflexive, and there is
an obvious induced map $\Quot^{\soft}(\ms E,\ell)\to\Tw^{ss}_{\ms
  X}(r,L,\Delta)$ whose image contains $p^{-1}([\ms E])$.  Using
Proposition \ref{sec:soft-quot-schemes-2}, we conclude that
$\dim p^{-1}([\ms E])\leq\ell(r+1)$.  If this is an equality, then we
see that a general soft quotient of $\ms E$ appears in $Y$.

\begin{lem}
  The kernel $\ms F$ of a general quotient $\ms E\to\ms Q$ satisfies
$$\dim\hom(\ms F,\ms F\tensor\omega_X)_0<\dim\hom(\ms E,\ms E\tensor\omega_X)_0.$$
\end{lem}
\begin{proof}
  First, we note that there is always a natural inclusion $\hom(\ms
  F,\ms F\tensor\omega_X)_0\inj\hom(\ms E,\ms E\tensor\omega_X)_0$
  arising from the fact that $\ms E$ is the soft hull of $\ms F$.

  Given $f\in\hom(\ms E,\ms E\tensor\omega_X)_0$, the fact that it is
  traceless shows that there is a collection of $\ell$ points $p_i$ of
  $X^{\sp}$ such that the fiber of $f$ at each $p_i$ is not a multiple
  of the identity (with respect to any local trivialization of
  $\omega_X$).  Thus, there is a quotient $\ms E(p_i)\to\ms Q_i$ of
  length $1$ whose kernel is not preserved by $f$.  The point $\ms
  E\to\oplus\ms Q_i$ thus has a kernel $\ms F$ such that $\hom(\ms
  F,\ms F\tensor\omega_X)_0$ does not contain $f$, as desired.
\end{proof}

Using the preceding lemma and Lemma \ref{L:bdry-codim}, we thus find
inequalities
$$\dim Z_1\geq\dim Y-\ell(r+1)\geq\dim Z-(r-1)-\ell(r+1)\geq\dim
Z-(2r-1)\ell-1,$$
and, moreover, if equality holds we know that
$\beta(Z)<\beta(Z_1)$. (In general, since $\beta$ is upper
semicontinuous, we know that $\beta(Z)\leq\beta(Z_1)$.)

Following Langer, write $\eps_i=\dim Z_i-(\dim Z_{i-1}-(2r-1)\ell-1)$,
so that $\beta(Z_{i-1})\leq\beta(Z_i)\leq\beta_\infty$ for all $i$,
and whenever $\eps_i=0$ the inequality is strict.  We conclude that
$\sum\eps_i\geq t-\beta_\infty$ and that 
$$\dim Z_t=\dim Z-(2r-1)\sum\ell_i-t+\sum\eps_i.$$
Moreover, since we are taking soft hulls, it is easy to see that
$\Delta=\Delta_t+2r\sum\ell_i$, from which we deduce that
$$\dim Z_t\geq\dim
Z-\left(1-\frac{1}{2r}\right)(\Delta-\Delta_t)-\beta_\infty.$$
Using the assumption about $\dim Z$, we thus conclude that
$$\dim Z_t-\left(1-\frac{1}{2r}\right)\Delta_t\geq
\left(\frac{1}{2r}-\frac{r-1}{2\theta(r)}\right)\Delta+C_1\sqrt{\Delta}+C_3-\beta_\infty.$$
Since ${\softboundary}Z_t^s=\emptyset$, we have by Theorem
\ref{T:bdry-exists} that
$$\dim
Z_t-\left(1-\frac{1}{2r}\right)\Delta_t<\left(\frac{1}{2r}-\frac{r-1}{2\theta(r)}\right)\Delta+C_1\sqrt{\Delta}+C_2.$$
Using the fact that $\Delta>\Delta_t$ and the description of $C_3$, we
immediately get a contradiction by comparing the two inequalities.
\end{proof}

\begin{cor}\label{sec:asympt-prop-moduli-2}
  There exists a constant $A_2\geq A_1$ such that whenever $\Delta\geq
  A_2$,
  \begin{enumerate}
  \item every irreducible component of $\Tw^{ss}_{\ms X}(r,L,\Delta)$
    is generically good, hence contains $\Tw^s_{\ms X}(r,L,\Delta)$ as
    a dense open substack which is generically smooth of the expected
    dimension;
  \item the stack $\Tw^{ss}_{\ms X}(r,L,\Delta)$ is normal is a
    normal local complete intersection.
  \end{enumerate}
\end{cor}
By ``local complete intersection'' we mean that for any smooth cover
by a scheme $U\to\Tw^{ss}_{\ms X}(r,L,\Delta)$, we have that $U$ is a
local complete intersection over the base.  (It is fruitful to compare
this with the classical statement -- as in Theorem 9.3.3 of \cite{h-l}
-- where the use of GIT means that only the stable locus is understood
to be lci, as the local structure at boundary points deteriorates in
the passage to the GIT quotient.)
\begin{proof}
  The first part is a direct result of Theorem \ref{T:good-big}.  To
  prove the second part, one could argue as in the proof of Theorem
  9.3.3 (and thus Theorem 2.2.8) of \cite{h-l} and invoke a smooth
  cover by a Quot scheme.  Instead, let us give an intrinsic proof.
  
  Let $\ms F$ be a semistable $\ms X$-twisted sheaf of determinant $L$
  and discriminant $\Delta$, and let $(V,v)\to\Tw^{ss}_{\ms
    X}(r,L,\Delta)$ be a smooth map from a pointed scheme whose
  completion at $v$ gives a miniversal deformation of $\ms F$ over $k$.  Since
  the good locus is dense, there will be points $w\in V(k)$ parametrizing
  good (hence stable) $\ms X$-twisted sheaves $\ms G$.  The expected
  dimension of $\Tw^s_{\ms X}(r,L,\Delta)$ at $[\ms G]$ is $\chi(\ms
  O)-\chi(\ms G,\ms G)$, while the expected dimension at $\ms F$ is
  $\chi(\ms O)-\chi(\ms F,\ms F)+\dim\Hom_0(\ms F,\ms F)$; the latter
  term is recognizable as the dimension of the tangent space to
  $\Aut_0(\ms F)$ (where this means automorphisms in the moduli
  problem, i.e., autmorphisms which preserve the determinant).

  By the assumption that $\ms G$ is good, we have that the dimension
  of $\Tw^{ss}_{\ms X}(r,L,\Delta)$ at $[\ms G]$ is equal to $\chi(\ms
  O)-\chi(\ms G,\ms G)$, so that $\dim_wV=\chi(\ms O)-\chi(\ms G,\ms
  G)+\dim(V/\Tw^{ss}_{\ms X}(r,L,\Delta))$, where the last term is the
  (constant) relative dimension of the smooth morphism
  $V\to\Tw^{ss}_{\ms X}(r,L,\Delta)$.  Since $\ms F$ and $\ms G$
  belong the same flat family, we see that $\chi(\ms F,\ms F)=\chi(\ms
  G,\ms G)$, from which we conclude that $\dim(V/\Tw^{ss}_{\ms
    X}(r,L,\Delta))=\dim\Hom_0(\ms F,\ms F)$.  Since
  $\dim_vV=\dim_wV$, we see that $V$ has the expected dimension at
  $v$.  Schlessinger's criterion (e.g., Proposition 2.A.11 of
  \cite{h-l}) implies that $V$ is lci at $v$, which implies (by
  definition) that $\Tw^{ss}_{\ms X}(r,L,\Delta)$ is lci at $\ms F$,
  as desired.

  Now that we know $\Tw^{ss}_{\ms X}(r,L,\Delta)$ is lci with singular
  locus of arbirarily large codimension (for large enough $\Delta$),
  we see that it is R1 and S2, and is thus normal by Serre's
  criterion.
\end{proof}

\subsection{Irreducibility}
\label{sec:asumpt-irred}

We conclude with two irreducibility results.  The first one is
significantly stronger, but only applies to the classical situation in
which $X\to X^\sp$ is an isomorphism.  The second uses asymptotic
normality to ensure that there is a geometrically irreducible
component of the stack for sufficiently large $\Delta$ and works
without any hypothesis on the stackiness of $X$; it uses a trick due
to de Jong and Starr and shown to the author by de Jong.

\begin{thm}\label{T:bdry-mandatory}
  There exists $A_3\geq A_2$ such that for any $\Delta\geq A_3$, any
  irreducible component of $\Tw^s_{\ms X/k}(r,L,\Delta)$ contains
  \begin{enumerate}
  \item a point $[\ms F]$ which is good and locally free, and
  \item a point $[\ms F]$ which is good such that $\ell^\soft(\ms F)=1$.
  \end{enumerate}
\end{thm}
\begin{proof} The standard proof of this fact (as in section 9.6 of
  \cite{h-l}) works.  We recall the numerical proof of the second, and
  leave the similar proof of the first to the reference (where it is
  spelled out in greater detail).

  To prove the second, let
  $A_3=A_2+2r\left(\frac{\beta_\infty}{r-1}+1\right)$.  If $\Delta\geq
  A_3$ then $\Delta\geq A_2$, and thus
  ${\softboundary}Z^s\neq\emptyset$.  Let $Y$ be an irreducible
  component of the dominant component of $\softboundary Z^s$; suppose
  the soft colength of the sheaves parametrized by $Y$ is $\ell$, and
  let $\rho:Y\to\Tw^{s}_{\ms X}(r,L,\Delta-2r\ell)$ be the morphism
  coming from formation of the soft hull.  Let $Z'$ be the closure of
  $\rho(Y)$.  Proposition \ref{sec:soft-quot-schemes-2} shows that the
  fibers of $\rho$ have dimension at most $\ell(r+1)$, so that $\dim
  Z'\geq\dim Y-\ell(r+1)$.  In addition, Lemma \ref{L:bdry-codim}
  shows that $\dim Y\geq\dim Z-(r-1)$.  Putting these together yields
$$\dim Z'\geq \exp\dim\Tw^{s}_{\ms
  X}(r,L,\Delta-2r\ell)+(\ell-1)(r-1),$$ since $\exp\dim\Tw^{ss}_{\ms
  X}(r,L,\Delta-2r\ell)=\dim\Tw^{s}_{\ms X}(r,L,\Delta)-2r\ell$ (since
the latter has the expected dimension).  Now, if $\Delta-2r\ell\geq
A_2$ then $Z'$ is generically good by Theorem \ref{T:good-big}, in
which case $\dim Z'\leq\exp\dim\Tw^{s}_{\ms X}(r,L,\Delta-2r\ell)$,
which shows that $(r-1)(\ell-1)=0$, so that $\ell=1$.  On the other
hand, if $\Delta-2r\ell<A_2$ then $2r\ell>\Delta-A_2\geq
A_3-A_2=2r\left(\frac{\beta_\infty}{r-1}+1\right)$, so that
$(r-1)(\ell-1)>\beta_\infty$ and $\dim Z'\geq\exp\dim\Tw^{s}_{\ms
  X}(r,L,\Delta-2r\ell)+\beta_\infty$, which contradicts Corollary
\ref{sec:dimension-estimate-7}.
\end{proof}

\begin{thm}\label{T:irred} The following asymptotic irreducibility
  statements hold.
  \begin{enumerate}
  \item If $X\to X^{\sp}$ is an isomorphism, then for all sufficiently
    large values of $\Delta$ the stack $\Tw^s_{\ms X/k}(r,L,\Delta)$
    is geometrically integral whenever it is non-empty, and if it is
    non-empty then so is $\Tw^s_{\ms X/k}(r,L,\Delta+1)$.
  \item In general, if there exists a totally regular locally free
    $\ms X$-twisted sheaf of rank $r$, then for any positive integer
    $E$, there is $\Delta>E$ such that $\Tw^s_{\ms X/k}(r,L,\Delta)$
    contains a geometrically integral connected component when it is
    non-empty.  Moreover, if it is non-empty for one value of $\Delta$
    then it is non-empty for infinitely many.
  \end{enumerate}
\end{thm}

\begin{proof}[Proof of Theorem \ref{T:irred}(1)]
  This works precisely as in the proof of Theorem 3.2.4.11 of
  \cite{twisted-moduli}, using the fact that the sheaves are totally
  regular.  We sketch the main points of the proof using a sequence of
  Propositions.

  \begin{prop}\label{sec:asympt-prop-moduli-3}
    Let $\ms E$ and $\ms F$ be two locally free totally regular $\ms
    X$-twisted sheaves of rank $r$ with determinant $L$.  For
    sufficiently large $m$, the cokernel
    of a general map $\ms E\to\ms F(m)$ is an invertible twisted sheaf
    supported on a smooth member of $|\ms O(rm)|$ in $X^\sp$.
  \end{prop}
  \begin{proof}
    We may choose $m$ sufficiently large that for any closed substack
    $Y\subset X$ of length at most $3$ the restriction map $\hom(\ms
    E,\ms F(m))\to\hom_Y(\ms E_Y,\ms F(m)_Y)$ is surjective.  Let $\A$
    be the affine space whose $k$-points are $\hom(\ms E,\ms F(m))$,
    and let $\Phi:\ms E_\A\to\ms F(m)_\A$ be the universal map over
    $\ms X\times\A$.  For a point $x\in X^\sp(k)$ write $Y_x=\spec\ms
    O_{X,x}/\mf m_x^2$.  There is a (non-canonical) isomorphism
    $\hom_{Y_x}(\ms E_{Y_x},\ms F(m)_{Y_x})\cong\M_r(\ms O_{X,x}/\mf
    m_x^2)$, where $r=\rk\ms E=\rk\ms F$.  We can write an element of
    the latter as $A=A_0+xA_1+yA_2$ with the $A_i$ elements of
    $\M_n(k)$.  We claim that the locus of $A$ such that $\det
    A=0\in\ms O_{X,x}/\mf m_x^2$ is a cone of codimension at least
    $3$.  To see this, we recall the Jacobi formula:
$$\det(A)=\det A_0+\operatorname{Tr}(\operatorname{adj}(A_0)(A_1x+A_2y)).$$
The condition that $\det A_0=0$ has codimension $1$.  It is clear from
the formula that the vanishing to first order in $x$ and $y$ will also
be conditions of codimension at least $1$.

Thus, there is a cone $U_x\subset\A$ of codimension at least $3$ such
that of $\alpha$ is not in $U_x$, then the map $\ms E\to\ms F(m))$
corresponding to $\alpha$ does not vanish in the second jet bundle at
$x$.  Since the $U_x$ vary algebraically with $x$, we see that there
is a locally closed subset $U\subset\A$ of codimension at least $1$ whose
complement $V$ parametrizes maps whose determinants are smooth at each
$x\in X^\sp(k)$.  The condition that the map be an isomorphism in the
fibers over $X\setminus X^\sp$ is open (and nonempty, since both $\ms
E$ and $\ms F$ are totally regular).  Finally, a similar argument
shows that the condition that the cokernel of $\ms E\to\ms F(m)$ have
rank $1$ at each point of its support is also a non-empty open subset.
Putting these together yields the Proposition.
\end{proof}
Now let $\ms Q$ be the cokernel of a general map as in Proposition
\ref{sec:asympt-prop-moduli-3}.  Write $\ms Q=i_{\ast}\ms M$, where
$i:C\to X^\sp$ is the inclusion of the (smooth) support of $\ms Q$
(which lies in $|rmH|$) and $\ms M$ is an invertible twisted sheaf on $C$.
Since $\ms E$ and $\ms F$ have the
same determinant and discriminant, we know the degree of $\ms M$.  (To
see this, one could for example pull back to a finite flat cover
$Y\to\ms X$ to reduce to the case where $\ms X$ is a smooth projective
surface.  Now knowing the determinant and discriminant of $\ms E$ and
$\ms F$ tells us their Hilbert polynomials, whence we know the Hilbert
polynomial of $\ms M$, which tells us the degree.)

This leads us to the second major step of the proof of Theorem
\ref{T:irred}, which is collected in the following pair of Propositions.

\begin{prop}\label{sec:asympt-prop-moduli-4}
  Let $\ms C\to C\to S$ be a $\m_N$-gerbe on a proper smooth family of
  geometrically connected curves over a connected locally Noetherian
  base scheme $S$.  Given a rational number $d$, the relative twisted
  Picard space $\Pic_{\ms C/S}(d)$ of
  invertible twisted sheaves on the fibers of $\ms C$ with degree $d$
  is a proper smooth algebraic $S$-space with geometrically integral fibers.
  Thus, if $S$ is integral then so is $\Pic_{\ms C/S}$.
\end{prop}
\begin{proof}
  The properness of $\Pic_{\ms C/S}(d)$ is a standard result, and
  works just as in the classical case (see
  section 3.1 of \cite{twisted-moduli} for details).  Thus, it remains
  to see that the geometric fibers are integral.  If
  $S=\spec\widebar\kappa$ then there is an invertible $\ms C$-twisted
  sheaf, say $\ms L$.  Twisting down by $\ms L$ gives an isomorphism
  between $\Pic_{\ms C/S}(d)$ and $\Pic_{C/S}(d')$ for a suitable
  integer $d'$.  Since the latter is a torsor under the Jacobian of
  $C$, we immediately see that it is integral, as desired.
\end{proof}

\begin{prop}\label{sec:asympt-prop-moduli-5}
  Let $\ms V$ and $\ms W$ be two locally free totally regular $\ms X$-twisted sheaves with the same
  rank, determinant, and discriminant.  There exists finite colength
  twisted subsheaves $\ms V'\subset\ms V$ and $\ms W'\subset\ms W$
  with the same colength and an irreducible family $\mf F$ of $\ms X$-twisted
  sheaves containing both $\ms V'$ and $\ms W'$.  Moreover, if $\ms V$
  and $\ms W$ are good, then we can assume that the members of $\mf F$
  are also good.
\end{prop}
\begin{proof}
  The proof is identical to the proof of Proposition 3.2.4.22 of
  \cite{twisted-moduli}.  The idea is as follows: choose $m$ so that
  there are exact sequences
$$0\to\ms V(-m)\to\ms W\to\ms P\to 0$$
and
$$0\to\ms V(-m)\to\ms V\to\ms Q\to 0,$$
where $\ms P$ and $\ms Q$ are invertible sheaves supported on smooth
members of $|rmH|$ lying in $X^\sp$.  Using Proposition
\ref{sec:asympt-prop-moduli-4}, we can connect $\ms P$ and $\ms Q$ in
an irreducible family, say $\ms D$ parametrized by $T$.

In general, given an extension $0\to\ms A\to\ms B\to\ms C\to 0$ with
$\ms C$ an invertible sheaf supported on a smooth curve, taking the
preimage of $\ms C(-\ell)$ under an inclusion $\ms C(-\ell)\to\ms C$
induced by a section of $\ms O(\ell)$ yields a finite colength
subsheaf $\ms B'\subset\ms B$.  Moreover, taking large values of
$\ell$ makes the space $\ext^1(\ms C(-\ell),\ms A)$ behave better in a
family (e.g., as $\ms C$ varies the $\ext^1$s form a vector bundle).
Taking $\ell$ sufficiently large that $\ext^1(\ms V(-m),\ms
D_t(-\ell))$ behaves well as $t$ varies, we get a family of sheaves interpolating
between $\ms V'$ and $\ms W'$, as desired.  (One goal of increasing
$\ell$ is to ensure that general such extensions are torsion free in
general fibers.  Once this is achieved, the good locus will be further
open subset, preserving irreducibility of the base of the family.)  We
refer the reader to the proof of [\emph{loc.\ cit\/}.] for details.
\end{proof}

Now the crux of the proof of Theorem \ref{T:irred} lies in Theorem
\ref{T:bdry-mandatory}.  (The mechanism is familiar from \cite{h-l},
\cite{langer-castelnuovo}, and \cite{twisted-moduli}.)  Let
$\Xi(r,L,\Delta)$ denote the set of irreducible components of
$\Tw^{ss}_{\ms X}(r,L,\Delta)$.  Suppose $\Delta\geq A_3$.  Given an
element $[Z]\in\Xi(r,L,\Delta)$, we then have that there is a point
$z\in Z$ parametrizing a good locally free $\ms X$-twisted sheaf $\ms
V$.  Choose a closed point $x\in X^\sp$ and let $\ms V'\subset\ms V$
be the kernel of a quotient $\ms V\to\ms Q$, where $\ms Q$ has length
$1$ and is supported at $x$.  The sheaf $\ms V'$ lies in a
well-defined component of $\Tw^{ss}_{\ms X}(r,L,\Delta+2r)$ because it
is good (and thus cannot lie in an intersection of components),
yielding a map $\phi(\Delta):\Xi(r,L,\Delta)\to\Xi(r,L,\Delta+2r)$.
Taking a further quotient of $\ms F$ at such a skyscraper twisted
sheaf produces an element $q\in\Xi(r,L,\Delta+4r)$ which is easily
seen to be equal to $\phi(\Delta+2r)\circ\phi(\Delta)([Z])$ (and
similarly for higher compositions and further quotients).

On the other hand, any element $[W]$ of $\Xi(r,L,\Delta+2r)$ contains
a point $w$ parametrizing a good torsion free $\ms X$-twisted sheaf
$\ms F$ with soft colength $1$ (and which is locally free in
neighborhoods of the
non-trivial residual gerbes).  Letting $\ms V:=\ms F^{\vee\vee}$, we
see that $[W]$ is the image of the component containing $\ms V$ under
the map $\phi(\Delta)$.  This shows that $\phi(\Delta)$ is surjective
for all $\Delta\geq A_3$.  This yields a sequence of surjections
$\Xi(r,L,\Delta)\surj\Xi(r,L,\Delta+2r)\surj\cdots$ for each $\Delta$.
Since $\Delta$ has bounded denominator, there are finitely many such
sequences containing all possible values of $\Delta\geq A_3$.

By Proposition \ref{sec:asympt-prop-moduli-5}, starting with $\ms V$
and $\ms W$, locally free good sheaves corresponding to components
$A$ and $B$ described by $\Xi(r,L,\Delta)$, there exists some $\ell$ and good colength
$\ell$ subsheaves $\ms V'\subset\ms V$ and $\ms W'\subset\ms W$ which
lie in the same irreducible component of $\Tw^{ss}_{\ms
  X}(r,L,\Delta+2r\ell)$.  This shows that the composition
$\phi(\Delta+2r\ell)\circ\phi(\Delta+2r(\ell-1))\circ\cdots\circ\phi(\Delta)$
collapses the two points $A$ and $B$ to the same point.  Since each
$\Xi(r,L,\Delta)$ is finite, we see that the sequence of surjections
eventually results in singletons.  Since there are only finitely many
such sequences to consider, we see that for $\Delta$ sufficiently
large, $\Xi(r,L,\Delta)$ is a singleton, i.e., that $\Tw^{ss}_{\ms
  X}(r,L,\Delta)$ is irreducible.
\end{proof}

\begin{proof}[Proof of Theorem \ref{T:irred}(2)] 
  We sketch the proof, leaving a few details to the reader.  A more
  complete treatment of a similar statement will appear in a
  forthcoming paper on the period-index problem for Brauer groups of
  surfaces over finite fields.

  Since extension to the perfect closure of $k$ is radicial, in order
  to prove the statement we may assume that the base field $k$ is
  perfect.  The absolute Galois group of $k$ acts on the components of
  $\Tw^s_{\ms X/k}(r,L,\Delta)\tensor_k\widebar k$ (for any chosen
  algebraic closure of $k$), and a component which is Galois-fixed
  corresponds to a geometrically integral component over $k$.

  For sufficiently large $\Delta$, the irreducible components of
  $\Tw^s_{\ms X/k}(r,L,\Delta)\tensor_k\widebar k$ are normal.
  Moreover, the map
  $\phi(\Delta):\Xi(r,L,\Delta)\to\Xi(r,L,\Delta+2r)$ defined above in
  the proof of Theorem \ref{T:irred}(1) is Galois equivariant.  Since
  $\Xi(r,L,\Delta)$ is finite, if we can show that iteration of the
  maps $\phi(\Delta)$ eventually brings any two points together, we
  will thus find a point in $\Xi(r,L,\Delta+2r\ell)$ (for some $\ell$)
  which is Galois fixed, yielding the desired component.  Since we no
  longer know that $\phi(\Delta)$ is surjective, we cannot conclude
  that this is the only component.

  The proof that the $\phi(\Delta)$ are contracting is similar to
  Proposition \ref{sec:asympt-prop-moduli-3}.  Given good torsion free
  $\ms X$-twisted sheaves $\ms E$ and $\ms F$, there are cofinite
  length subsheaves $\ms E'$ and $\ms F'$ which have isomorphic fibers
  over every Henselization of $X^{\sp}$.  Now the analogue of Proposition
  \ref{sec:asympt-prop-moduli-3} for torsion free $\ms X$-twisted
  sheaves yields further cofinite length subsheaves $\ms E''\subset\ms
  E'$ and $\ms F''\subset\ms F'$ which lie in a connected family.
  Since $\ms E''$ and $\ms F''$ are good, this shows that the
  components containing $\ms E'$ and $\ms F'$ are contracted by a
  suitable iterate of the maps $\phi(\Delta)$.
\end{proof}

\end{document}